\documentclass[12pt]{article}

\usepackage[english]{babel}
\usepackage[T1]{fontenc}
\usepackage[applemac]{inputenc}
\usepackage{lmodern}
\usepackage{microtype}
\usepackage{bbm}
\usepackage{tikz}
\usepackage{graphicx}
\usepackage{subfloat}
\usepackage{amsmath}
\usepackage{amssymb}
\usepackage{amsthm}
\usepackage{latexsym}
\usepackage{color}
\usepackage{dsfont}
\usepackage{appendix}

\usepackage{hyperref}
\usepackage{float}
\usepackage{enumitem}
\usepackage{mathabx}
\usepackage{wrapfig}
\usepackage{graphicx} 
\graphicspath{{figures/}}

\DeclareSymbolFont{calletters}{OMS}{cmsy}{m}{n}
\DeclareSymbolFontAlphabet{\mathcal}{calletters}

\def\be{\begin{eqnarray}}
\def\ee{\end{eqnarray}}

\def\b*{\begin{eqnarray*}}
\def\e*{\end{eqnarray*}}

\def \C{\mathbb{C}}

\def \L{\mathbb{L}}

\def \N{\mathbb{N}}
\def \P{{\mathbb P}}

\def \R{\mathbb{R}}

\def \a{\alpha}

\def \eps{\varepsilon}

\def\Ac{{\cal A}}
\def\Bc{{\cal B}}
\def\Cc{{\cal C}}
\def\Dc{{\cal D}}
\def\Ec{{\cal E}}

\def\Kc{{\cal K}}

\def\Mc{{\cal M}}
\def\Nc{{\cal N}}

\def\Pc{{\cal P}}
\def\Qc{{\cal Q}}

\def\Sc{{\cal S}}
\def\Tc{{\cal T}}

\def\Xc{{\cal X}}
\def\Yc{{\cal Y}}

\def\Bb{\overline{B}}

\def\Ib{\overline{I}}

\def \ri{{\rm ri\hspace{0cm}}}

\def \cl{{\rm cl\hspace{0.05cm}}}
\def \interior{{\rm int\hspace{0cm}}}
\def \dom{{\rm dom}}
\def \supp{{\rm supp}}

\def \conv{{\rm conv}}

\def \Aff{{\rm Aff}}

\def \aff{{\rm aff}}

\def \Leb{{\L}}
\def \Ctn{{\rm{C}}}

\def \Ibf{{\mathbf I}}
\def \Sbf{{\mathbf S}}

\def \rank{{\rm rank}}

\def \Kcirc{{\wideparen{\Kc}}}

\def\Ima{{\rm Im}}

\def\no{\noindent}
\def\x{\times}

\def\={\;=\;}
\def\.{\;.}

\def\eps{\varepsilon}

\def \1{{\bf 1}}
\def \ep{\hbox{ }\hfill{ ${\cal t}$~\hspace{-5.1mm}~${\cal u}$   } }

\def \proof{{\noindent \bf Proof. }}
\def \ep{\hbox{ }\hfill$\Box$}

 \def\normeL2#1{\left\|{#1}\right\|_{L^2}}

\title{Local structure of multi-dimensional martingale optimal transport\thanks{The author gratefully acknowledges the financial support of the ERC 321111 Rofirm, and the Chairs Financial Risks (Risk Foundation, sponsored by Soci\'et\'e G\'en\'erale) and Finance and Sustainable Development (IEF sponsored by EDF and CA).}}

\author{Hadrien De March\thanks{CMAP, \'Ecole Polytechnique, hadrien.de-march@polytechnique.org.}}

\date{\today}

\begin{document}

\maketitle

\newtheorem{Theorem}{Theorem}[section]
\newtheorem{Lemma}[Theorem]{Lemma}
\newtheorem{Corollary}[Theorem]{Corollary}
\newtheorem{Proposition}[Theorem]{Proposition}
\newtheorem{Remark}[Theorem]{Remark}
\newtheorem{Example}[Theorem]{Example}
\newtheorem{Definition}[Theorem]{Definition}
\newtheorem{Assumption}[Theorem]{Assumption}


\abstract{This paper analyzes the support of the conditional distribution of optimal martingale transport couplings between marginals in $\R^d$ for arbitrary dimension $d\ge 1$. In the context of a distance cost in dimension larger than $2$, previous results established by Ghoussoub, Kim \& Lim \cite{ghoussoub2015structure} show that this conditional distribution is concentrated on its own Choquet boundary. Moreover, when the target measure is atomic, they prove that the support of this distribution is concentrated on $d+1$ points, and conjecture that this result is valid for arbitrary target measure. 

We provide a structure result of the support of the conditional distribution for general Lipschitz costs. Using tools from algebraic geometry, we provide sufficient conditions for finiteness of this conditional support, together with (optimal) lower bounds on the maximal cardinality for a given cost function. More results are obtained for specific examples of cost functions based on distance functions. In particular, we show that the above conjecture of Ghoussoub, Kim \& Lim is not valid beyond the context of atomic target distributions. 

\vspace{1cm}

\noindent {\bf Key words.}  Martingale optimal transport, local structure, differential structure, support.
}

\section{Introduction}

The problem of martingale optimal transport was introduced as the dual of the problem of robust (model-free) superhedging of exotic derivatives in financial mathematics, see Beiglb\"ock, Henry-Labord\`ere \& Penkner \cite{beiglbock2013model} in discrete time, and Galichon, Henry-Labord\`ere \& Touzi \cite{galichon2014stochastic} in continuous-time. Previously the robust superhedging problem was introduced  by Hobson \cite{hobson1998robust}, and was addressing specific examples of exotic derivatives by means of corresponding solutions of the Skorokhod embedding problem, see \cite{cox2011robust,hobson2015robust,hobson2012robust}, and the survey \cite{hobson2011skorokhod}. 

Our interest in the present paper is on the multi-dimensional martingale optimal transport. Given two probability measures $\mu,\nu$ on $\R^d$, with finite first order moment, martingale optimal transport differs from standard optimal transport in that the set of all interpolating probability measures $\Pc(\mu,\nu)$ on the product space is reduced to the subset $\Mc(\mu,\nu)$ restricted by the martingale condition. We recall from Strassen \cite{strassen1965existence} that $\Mc(\mu,\nu)\neq\emptyset$ if and only if $\mu\preceq\nu$ in the convex order, i.e. $\mu(f)\le\nu(f)$ for all convex functions $f$. Notice that the inequality $\mu(f)\le\nu(f)$ is a direct consequence of the Jensen inequality, the reverse implication follows from the Hahn-Banach theorem.

This paper focuses on showing the differential structure of the support of optimal probabilities for the martingale optimal transport Problem. In the case of optimal transport, a classical result by R\"uschendorf \cite{ruschendorf1991frechet} states that if the map $y\longmapsto c_x(x_0,y)$ is injective, then the optimal transport is unique and supported on a graph, i.e. we may find $T:\Xc\longrightarrow\Yc$ such that $\P^*[Y=T(X)]=1$ for all optimal coupling $\P^*\in\Pc(\mu,\nu)$. The corresponding result in the context of the one-dimensional martingale transport problem was obtained by Beiglb\"ock-Juillet \cite{beiglboeck2016problem}, and further extended by Henry-Labord\`ere \& Touzi \cite{henry2016explicit}. Namely, under the so-called martingale Spence-Mirrlees condition, $c_x$ strictly convex in $y$, the left-curtain transport plan is optimal and concentrated on two graphs, i.e. we may find $T_d,T_u:\Xc\longrightarrow\Yc$ such that $\P^*[Y\in\{T_d(X),T_u(X)\}]=1$ for all optimal coupling $\P^*\in\Mc(\mu,\nu)$. In this case we get similarly the uniqueness by a convexity argument.

An important issue in optimal transport is the existence and the characterization of optimal transport maps. Under the so-called twist condition (also called Spence-Mirrlees condition in the economics litterature) it was proved that the optimal transport is supported on one graph. In the context of martingale optimal transport on the line, Beiglb\"ock \& Juillet introduced the left-monotone martingale interpolating measure as a remarkable transport plan supported on two graphs, and prove its optimality for some classes of cost functions. Ghoussoub, Kim \& Lim conjectured that in higher dimensional Martingale Optimal Transport for distance cost, the optimal plans will be supported on $d+1$ graphs. We prove here that there is no hope of extending this property beyond the case of atomic measure. This is obtained using the reciprocal property of the structure theorem of this paper, which serves as a counterexample generator. We further prove that for "almost all" smooth cost function, the optimal coupling are always concentrated on a finite number of graphs, and we may always find densities $\mu$ and $\nu$ that are dominated by the Lebesgue measure such that the optimal coupling is concentrated on $d+2$ maps for $d$ even.

A first such study in higher dimension was performed by Lim \cite{lim2014optimal} under radial symmetry that allows in fact to reduce the problem to one-dimension. A more "higher-dimensional specific" approach was achieved by Ghoussoub, Kim \& Lim \cite{ghoussoub2015structure}. Their main structure result is that for the Euclidean distance cost, the supports of optimal kernels will be concentrated on their own Choquet boundary (i.e. the extreme points of the closure of their convex hull).

Our subsequent results differ from \cite{ghoussoub2015structure} from two perspectives. First, we prove that with the same techniques we can easily prove much more precise results on the local structure of the optimal Kernel, in particular, we prove that they are concentrated on $2d$ (possibly degenerate) graphs, which is much more precise than a concentration on the Choquet boundary. Our main structure result states that the optimal kernels are supported on the intersection of the graph of the partial gradient $c_x(x_0,\cdot)$ with the graph of an affine function $A_{x_0}\in\Aff_d$. Second, we prove a reciprocal property, i.e. that for any subset of such graph intersection $\{c_x(x_0,Y)=A(Y)\}$ for $A\in\Aff_d$, we may find marginals such that this set is an optimizer for these marginals. Thanks to this reciprocal property we prove that Conjecture 2 in \cite{ghoussoub2015structure} that we mentioned above is wrong. They prove this conjecture in the particular case in which the second marginal $\nu$ is atomic, however in view of our results it only works in this particular case, as we produce counterexamples in which $\mu$ and $\nu$ are dominated by the Lebesgue measure. Indeed, we prove that the support of the conditional kernel is characterized by an algebraic structure independent from the support of $\nu$, then when this support is atomic, very particular phenomena happen. Thus the intuition suggests that finding this kind of solution for an atomic approximation of a non-atomic $\nu$ is not a stable approach, as in the limit there are generally $2d$ points in the kernel.

The paper is organized as follows. Section \ref{sect:main} gives the main results: Subsection \ref{subsect:structure} states the Assumption and the main structure theorem, Subsection \ref{subsect:geoalg} applies this theorem to show the relation between finiteness of the conditional support and the algebraic geometry of its derivatives, Subsection \ref{subsect:lowerbound} gives the maximal cardinality that is universally reachable for the support up to choosing carefully the marginals, and finally Subsection \ref{subsect:charactclassical} shows how the structure theorem applied to classical costs like powers of the Euclidean distance allows to give precise descriptions and properties of the conditional supports of optimal plans. Finally Section \ref{sect:proof} contains all the proofs to the results in the previous sections, and Section \ref{sect:numerics} provides some numerical experiments.
\\

\no {\bf Notation}\quad We fix an integer $d\ge 1$. For $x\in\R$, we denote $sg(x):=\mathbf{1}_{x>0}-\mathbf{1}_{x<0}$. If $f:\R\longrightarrow\R$ we denote by $\mathbf{fix}(f)$ the set of fixed points of $f$. A function $f:\R^d\longrightarrow \R^d$ is said to be super-linear if $\lim_{|y|\to\infty}\frac{|f(y)|}{|y|} = \infty$. Let a function $f:\R^d\longrightarrow \R$ and $x_0\in\R^d$, we say that $f$ is super-differentiable (resp. sub-differentiable) at $x_0$ if we may find $p\in\R^d$ such that $f(x)-f(x_0)\le p\cdot(x-x_0) + o(x-x_0)$ (resp. $\ge$) when $x\longmapsto x_0$, in this condition, we say that $p$ belongs to the super-gradient $\partial^+f(x_0)$ (resp. sub-gradient $\partial^-f(x_0)$) of $f$ at $x_0$. This local notion extends the classical global notion of super-differential (resp. sub) for concave (resp. convex) functions.

For $x\in\R^d$, $r\ge 0$, and $V$ an affine subspace of dimension $d'$ containing $x$, we denote $\Sc_V(x,r)$ the $\dim V-1$ dimensional sphere in the affine space $V$ for the Euclidean distance, centered in $x$ with radius $r$. We denote by $\Aff_d$ the set of Affine maps from $\R^d$ to itself. Let $A\in\Aff_d$, notice that its derivative $\nabla A$ is constant over $\R^d$, we abuse notation and denote $\nabla A$ for the matrix representation of this derivative. Let $M\in\Mc_d(\R)$, a real matrix of size $d$, we denote $\det M$ the determinant of $M$, $\ker M$ is the kernel of $M$, $\Ima M$ is the image of this matrix, and $Sp(M)$ is the set of all complex eigenvalues of $M$. We also denote $Com(M)$ the comatrix of $M$: for $1\le i,j\le d$, $Com(M)_{i,j} = (-1)^{i+j}\det M^{i,j}$, where $M^{i,j}$ is the matrix of size $d-1$ obtained by removing the $i^{th}$ line and the $j^{th}$ row of $M$. Recall the useful comatrix formula:
\be\label{eq:comatrix}
Com(M)^t M = M Com(M)^t = (\det M) I_d.
\ee

As a consequence, whenever $M$ is invertible, $M^{-1} = \frac{1}{\det M}Com(M)^t$. Throughout this paper, $\R^d$ is endowed with the Euclidean structure, the Euclidean norm of $x\in\R^d$ will be denoted $|x|$, the $p-$norm of $x$ will be denoted $|x|_p:=\left(\sum_{i=1}^d|x_i|^p\right)^\frac{1}{p}$. We denote $(e_i)_{1\le i\le d}$ the canonical basis of $\R^d$. Let $B\subset E$ with $E$ a vector space, we denote $B^*:=B\setminus\{0\}$, and $|B|$ the possibly infinite cardinal of $B$. If $V$ is a topological affine space and $B\subset V$ is a subset of $V$, $\interior B$ is the interior of $B$, $\cl B$ is the closure of $B$, $\aff B$ is the smallest affine subspace of $V$ containing $B$, $\conv B$ is the convex hull of $B$, $\dim(B):=\dim(\aff B)$, and $\ri B$ is the relative interior of $B$, which is the interior of $B$ in the topology of $\aff B$ induced by the topology of $V$. We also denote by $\partial B:=\cl B\setminus\ri B$ the relative boundary of $B$, and if $V$ is endowed with a euclidean structure, we denote by $proj_B(x)$ the orthogonal projection of $x\in V$ on $\aff B$. A set $B$ is said to be discrete if it consists of isolated points.

We denote $\Omega:=\R^d\times\R^d$ and define the two canonical maps
\b*
 X :(x,y)\in\Omega
 \longmapsto x\in\R^d
 &\mbox{and}&
 Y :(x,y)\in\Omega
 \longmapsto y\in\R^d.
 \e*
For $\varphi,\psi:\R^d\longrightarrow\bar\R$, and $h:\R^d\longrightarrow\R^d$, we denote 
\b*
\varphi\oplus\psi
:=
\varphi(X)+\psi(Y),
&\mbox{and}&
h^\otimes := h(X)\cdot(Y-X),
\e*
with the convention $\infty-\infty = \infty$.

For a Polish space $\Xc$, we denote by $\Pc(\Xc)$ the set of all probability measures on $\big(\Xc,\Bc(\Xc)\big)$. For $\P\in\Pc(\Xc)$, we denote by $\supp\P$ the smallest closed support of $\P$. Let $\Yc$ be another Polish space, and $\P\in\Pc(\Xc\x\Yc)$. The corresponding conditional kernel $\P_x$ is defined by:
$$\P(dx,dy) = \mu(dx) \P_x(dy),\text{ where }\mu:=\P\circ X^{-1}.$$

Let $n\ge 0$ and a field $\mathbb{K}$ ($\R$ or $\C$ in this paper), we denote $\mathbb{K}_n[X]$ the collection of all polynomials on $\mathbb{K}$ of degree at most $n$. The set $\C^{hom}[X]$ is the collection of homogeneous polynomials of $\C[X]$. Similarly for $k\ge 1$, we define $\mathbb{K}_n[X_1,...,X_d]$ the collection of multivariate polynomials on $\mathbb{K}$ of degree at most $n$. We denote the monomial $X^\alpha := X_1^{\alpha_1}...X_d^{\alpha_d}$, and $|\alpha|=\alpha_1+...+\alpha_d$ for all integer vector $\alpha\in\N^d$. For two polynomial $P$ and $Q$, we denote $\gcd(P,Q)$ their greatest common divider. Finally, we denote $\P^d := \big(\C^{d+1}\big)^*/\C^*$ the projective plan of degree $d$.
\\

\no {\bf The martingale optimal transport problem}\quad Throughout this paper, we consider two probability measures $\mu$ and $\nu$ on $\R^d$ with finite first order moment, and $\mu \preceq \nu$ in the convex order, i.e. $\nu(f)\ge \mu(f)$ for all integrable convex $f$. We denote by $\Mc(\mu,\nu)$ the collection of all probability measures on $\R^d\times\R^d$ with marginals $\P\circ X^{-1}=\mu$ and $\P\circ Y^{-1}=\nu$. Notice that $\Mc(\mu,\nu)\neq\emptyset$ by Strassen \cite{strassen1965existence}.

An $\Mc(\mu,\nu)-$polar set is an element of $\cap_{\P\in\Mc(\mu,\nu)}\Nc_\P$. A property is said to hold $\Mc(\mu,\nu)-$quasi surely (abbreviated as q.s.) if it holds on the complement of an $\Mc(\mu,\nu)-$polar set.

For a derivative contract defined by a non-negative cost function $c:\R^d\times\R^d\longrightarrow\R_+$, the martingale optimal transport problem is defined by:
 \be\label{pb:MOT}
 \Sbf_{\mu,\nu}(c)
 &:=&
 \sup_{\P\in\Mc(\mu,\nu)}
 \P[c].
 \ee

The corresponding robust superhedging problem is
 \be\label{Robustsuperhedge}
 \Ibf_{\mu,\nu}(c)
 &:=&
 \inf_{(\varphi,\psi,h)\in\Dc_{\mu,\nu}(c)} \mu(\varphi)+\nu(\psi),
 \ee
where
 \be
 \Dc_{\mu,\nu}(c)
 &:=&
 \big\{ (\varphi,\psi,h)\in\L^1(\mu)\x\L^1(\nu)\x\L^1(\mu,\R^d):
                                                           ~\varphi\oplus\psi+h^\otimes\ge c
 \big\}.~~~~~
 \ee
 The following inequality is immediate:
  \be\label{eq:weakduality}
\Sbf_{\mu,\nu}(c) \le \Ibf_{\mu,\nu}(c).
\ee
This inequality is the so-called weak duality. For upper semi-continuous cost, Beiglb\"ock, Henry-Labord\`ere, and Penckner \cite{beiglbock2013model}, and Zaev \cite{zaev2015monge} proved that strong duality holds, i.e. $\Sbf_{\mu,\nu}(c)= \Ibf_{\mu,\nu}(c)$. For any Borel cost function, De March \cite{de2018quasi} extended the quasi sure duality result to the multi-dimensional context, and proved the existence of a dual minimizer.

\section{Main results}\label{sect:main}

\subsection{Main structure theorem}\label{subsect:structure}

An important question in optimal transport theory is the structure of the support of the conditional distribution of optimal transport plans. Theorem \ref{thm:structure} below gives a partial structure to this question. As a preparation we introduce a technical assumption.

We denote $\Kcirc$ the collection of closed convex subsets of $\R^d$, which is a Polish space when endowed with the Wijsman topology (see Beer \cite{beer1991polish}). De March \& Touzi \cite{de2017irreducible} proved that we may find a Borel mapping $I:\R^d\longmapsto\Kcirc$ such that $\{I(x):x\in\R^d\}$ is a partition of $\R^d$, $Y\in\cl I(X)$, $\Mc(\mu,\nu)-$a.s. and $\cl I(X) = \cl\,\conv\,\supp\hat\P_X$, $\mu-$a.s.
 for some $\hat\P\in\Mc(\mu,\nu)$. As the map $I$ is Borel, $I(X)$ is a random variable, let $\eta:=\mu\circ I^{-1}$ be the push forward of $\mu$ by $I$. It was proved in \cite{de2018quasi} that the optimal transport disintegrates on all the "components" $I(X)$. The following conditions are needed throughout this paper.
 
\begin{Assumption}\label{ass:duality}
\no{\rm (i)} $c:\Omega\longmapsto \R$ is upper semi-analytic, $\mu\preceq\nu$ in convex order in $\Pc(\R^d)$, $c\ge \alpha\oplus \beta+\gamma^\otimes$ for some $(\a,\beta,\gamma)\in\Leb^1(\mu)\times\Leb^1(\nu)\x\Leb^0(\R^d,\R^d)$, and $\Sbf_{\mu,\nu}(c)<\infty$.

\no{\rm (ii)} The cost $c$ is locally Lipschitz and sub-differentiable in the first variable $x\in I$, uniformly in the second variable $y\in \cl I$, $\eta-$a.s.

\no{\rm (iii)} The conditional probability $\mu_I:=\mu\circ \big(X|I\big)^{-1}$ is dominated by the Lebesgue measure on $I$, $\eta-$a.s.
\end{Assumption}

The statements (i) and (ii) of Assumption \ref{ass:duality} are verified for example if $c$ is differentiable and if $\mu$ and $\nu$ are compactly supported. On another hand, the statement (iii) is much more tricky. It is well known that Sudakov \cite{sudakov1979geometric} thought that he had solved the Monge optimal transport problem by using the (wrong) fact that the disintegration of the Lebesgue measure on a partition of convex sets would be dominated by the Lebesgue measure on each of these convex sets. However, \cite{ambrosio2004existence}, provides a counterexample inspired from another paradoxal counterexample by Davies \cite{davies1952accessibility}. This Nikodym set $N$ is equal to the tridimensional cube up to a Lebesgue negligible set. Furthermore it is designed so that a continuum of mutually disjoint lines which intersect all $N$ in one singleton each. Thus the Lebesgue measure on the cube disintegrates on this continuum of lines into Dirac measures on each lines.

Statement (iii) is implied for example by the domination of $\mu$ by the Lebesgue measure together with the fact that $\dim I(X)\in\{0,d-1,d\}$, $\mu-$a.s. (see Lemma C.1 of \cite{ghoussoub2015structure} implying that the Lebesgue measure disintegrates in measures dominated by Lebesgue on the $d-1-$dimensional components), in particular together with the fact that $d\le 2$, or together with the fact that $\nu$ is the law of $X_\tau:=X_0+\int_0^t \sigma_sdW_s$, where $X_0\sim \mu$, $W$ a $d-$dimensional Brownian motion independent of $X_0$, $\tau$ is a positive bounded stopping time, and $(\sigma_t)_{t\ge 0}$ is a bounded cadlag process with values in $\Mc_{d}(\R)$ adapted to the $W-$filtration with $\sigma_0$ invertible. See the proof of Remark 4.3 in \cite{de2018quasi}.

\begin{Theorem}\label{thm:structure}
\no{\rm (i)} Under Assumption \ref{ass:duality} we may find $(A_x)_{x\in\R^d}\subset \Aff_d$ such that for all $\P^*\in\Mc(\mu,\nu)$ optimal for \eqref{pb:MOT},
\b*
x\in\ri\,\conv\,\supp\,\P^*_x,&\mbox{and}\quad\supp\,\P^*_{x}\subset\{ c_x(x,Y)=A_{x}(Y)\}&\mbox{for }\mu-\mbox{a.e. }x\in\R^d.
\e*

\no{\rm (ii)} Conversely, let a compact $S_0\subset \{c_x(x_0,Y)=A(Y)\}$ for some $x_0\in\R^d$ and $A\in\Aff_d$, be such that $x_0\in\interior\,\conv \,S_0$, $c$ is $\Ctn^{2,0}\cap\Ctn^{1,1}$ in the neighborhood of $\{x_0\}\x S_0$, and $c_{xy}(S_0)-\nabla A\subset GL_d(\R)$, then $S_0$ has a finite cardinal $k\ge d+1$ and we may find $\mu_0,\nu_0\in\Pc(\R^d)$ with $\Ctn^1$ densities such that
$$\P^*(dx,dy):=\mu_0(dx)\sum_{i=1}^k\lambda_i(x)\delta_{T_i(x)}(dy)$$

is the unique solution to \eqref{pb:MOT}, with $(T_i)_{1\le i\le k}\subset\Ctn^1(\supp\,\mu_0,\R^d)$ such that $S_0=\{T_i(x_0)\}_{1\le i\le k}$, and $(\lambda_i)_{1\le i\le k}\subset\Ctn^1(\supp\,\mu_0)$.
\end{Theorem}


\begin{Remark}
We have $\nabla A_x(x)=\nabla\varphi(x)-h(x)$ in Theorem \ref{thm:structure} from its proof. Under the stronger assumption that $\varphi$ and $h$ are $\Ctn^1$, we can get this result much easier. As for $(x,y)\in\R^d$,
$$\varphi(x)+\psi(y)+h(x)\cdot(y-x)-c(x,y)\geq 0,$$
with equality for $(x,y)\in\Gamma$. When $y_0$ is fixed, $x_0$ such that $(x_0,y_0)\in\Gamma$ is a critical point of $x\mapsto\varphi(x)+\psi(y_0)+h(x)\cdot(y_0-x)-c(x,y_0)$. Then we get
$c_x(x_0,y_0)=\nabla h(x_0)(y_0-x_0) + \nabla\varphi(x_0)-h(x_0)$ by the first order condition.

We see that we have in this case $A_{x_0}(y):=\nabla h(x_0)(y-x_0) + \nabla\varphi(x_0)-h(x_0)$, and $\Gamma_{x_0}\subset \{c_x(x_0,Y)=A_{x_0}(Y)\}$, for $\mu-$a.e. $x_0\in\R^d$.
\end{Remark}

\begin{Remark}\label{rmk:d+1determine}
Even though the set $S_0:=\{c_x(x_0,Y) = A(Y)\}$ for $x_0\in \R^d$ and $A\in\Aff_d$ may contain more than $d+1$ points, it is completely determined by $d+1$ affine independent points $y_1,...,y_{d+1}\in S_0$, as the equations $c_x(x_0,y_i) = A(y_i)$ determine completely the affine map $A$.
\end{Remark}

\no{\bf Proof of Theorem \ref{thm:structure}}
\no\underline{(i)} By Theorem 3.5 (i) in \cite{de2018quasi}, (and using the notation therein), the quasi-sure robust super-hedging problem may be decomposed in pointwise robust super-hedging separate problems attached to each components, and we may find functions $(\varphi,h)\in\Leb^0(\R^d)\x\Leb^0(\R^d,\R^d)$, and $(\psi_K)_{K\in I(\R^d)}\subset\Leb_+^0(\R^d)$ with $\psi_{I(X)}(Y)\in\Leb^0_+(\Omega)$, and $\dom\,\psi_I = J_\theta$, $\eta-$a.s. for some $\theta\in\widehat{\Tc}(\mu,\nu)$, such that $c
 \le
\varphi(X)+\psi_{I(X)}(Y)+h^\otimes$, and $
\Sbf_{\mu,\nu}(c)
 =
\Sbf_{\mu,\nu}\big(\varphi(X)+\psi_{I(X)}(Y)+h^\otimes\big)$. Then applying the theorem to $c':=\varphi(X)+\psi_{I(X)}(Y)+h^\otimes$, $\Sbf_{\mu,\nu}(c) = \Sbf_{\mu,\nu}\big(\varphi(X)+\psi_{I(X)}(Y)+h^\otimes\big) = \sup_{\P\in\Mc(\mu,\nu)}\Sbf_{\mu_I,\nu_I^\P}\big(\varphi(X)+\psi_{I(X)}(Y)+h^\otimes\big).$ Then if $\P\in\Mc(\mu,\nu)$ is optimal for $\Sbf_{\mu,\nu}(c)$, then $\P_I[c = \varphi\oplus\psi_{I}+h^\otimes]=1$, $\eta-$a.s. By Lemma 3.17 in \cite{de2018quasi} the regularity of $c$ in Assumption \ref{ass:duality} (ii) guarantees that we may chose $\varphi$ to be locally Lipschitz on $I$, and $h$ locally bounded on $I$. In view of Assumption \ref{ass:duality} (iii), $\varphi$ is differentiable $\mu_I-$a.e. by the Rademacher Theorem. Then after possibly restricting to an irreducible component, we may suppose that we have the following duality: for any $x,y\in \R^d$,
\be\label{eq:pwdual}
\varphi(x)+\psi(y)+h(x)\cdot(y-x)-c(x,y) &\geq& 0,
\ee
with equality if and only if $(x,y)\in \Gamma:=\{\varphi\oplus\psi+h^\otimes=c<\infty\}$, concentrating all optimal coupling for $\Sbf_{\mu,\nu}(c)$.

Let $x_0\in\ri\,\conv\,\dom\,\psi$ such that $\varphi$ is differentiable in $x_0$. Let $y_1,...,y_k\in \Gamma_{x_0}$ such that $\sum_{i=1}^k\lambda_i y_i = x_0$, convex combination. We complete $(y_1,...,y_k)$ in a barycentric basis $(y_1,...,y_k,y_{k+1},...,y_l)$ of $\ri\,\conv\,\dom\,\psi$. Let $x\in\ri\,\conv\,\dom\,\psi$ in the neighborhood of $x_0$, and let $(\lambda_i')$ such that $x=\sum_{i=1}^l\lambda_i' y_i $, convex combination. We apply \eqref{eq:pwdual}, both in the equality and in the inequality case:
\b*
\varphi(x)+\sum_{i=1}^l\lambda_i'\psi(y_i) \geq \sum_{i=1}^l\lambda_i'c(x,y_i),&\varphi(x_0)+\sum_{i=1}^l\lambda_i'\psi(y_i) + h(x_0)\cdot(x-x_0) = \sum_{i=1}^l\lambda_i'c(x_0,y_i).
\e*

By subtracting these equations, we get
\b*
\varphi(x)-\varphi(x_0)-h(x_0)\cdot(x-x_0)&\ge& \sum_{i=1}^l\lambda_i'\big(c(x,y_i)-c(x_0,y_i)\big).
\e*
As $c$ is Lipschitz in $x$, and $\lambda_i'\longrightarrow \lambda_i$ when $x\to x_0$, we get:
\b*
\big(\nabla\varphi(x_0)-h(x_0)\big)\cdot(x-x_0)+o(x-x_0)\ge \sum_{i=1}^k\lambda_i\big(c(x,y_i)-c(x_0,y_i)\big).
\e*
Then, $x\longmapsto \sum_{i=1}^k\lambda_ic(x,y_i)$ is super-differentiable at $x_0$, and $\nabla\varphi(x_0)-h(x_0)$ belongs to its super-gradient. As $x\longmapsto c(x,y)$ is sub-differentiable by Assumption \ref{ass:duality} (ii), it implies that $x\longmapsto c(x,y_i)$ is differentiable at $x_0$ for all $i$ such that $\lambda_i>0$, and therefore
\be\label{eq:sumdiff}
\nabla\varphi(x_0)-h(x_0)&=& \sum_{i=1}^k\lambda_ic_x(x_0,y_i).
\ee
Now we want to prove that we may find $A_x\in\Aff_d$ such that $A_x(y)=c_x(x,y)$ for all $y\in\Gamma_x$.

Let $y_1^0,...,y_m^0\in\Gamma_{x_0}$ generating $\aff\Gamma_{x_0}$ and such that $x\in\ri\,\conv(y_1^0,...,y_m^0)$, let $y\in \Gamma_{x_0}$. $A_x$ is defined in a unique way if $\nabla A = 0$ on $(\aff \Gamma_{x_0}-x_0)^\perp$ by its values on $(y_1^0,...,y_m^0)$. Now we prove that $A_x(y) = c_x(x_0,y)$. As $y\in\aff(y_1^0,...,y_m^0)$, we may find $(\mu_i)$ so that $\sum_{i=1}\mu_i y_i^0 = y$, and $\sum_{i=1}\mu_i = 1$. For $\eps>0$ small enough, $x_0-\eps(y-x_0)\in\ri\,\conv(y_1^0,...,y_m^0)$. Then $x_0-\eps(y-x_0) = \sum_{i=1}\lambda_i^0 y_i$ with $\lambda_i^0>0$. We take the convex combination: $x_0=\frac{1}{1+\eps}(x_0-\eps(y-x_0))+\frac{\eps}{1+\eps}y$, and $x_0=\sum_{i=1}\big(\frac{1}{1+\eps}\lambda_i^0+\frac{\eps}{1+\eps}\mu_i\big) y_i^0$. We suppose that $\eps$ is small enough so that $\lambda^\eps_i:=\frac{1}{1+\eps}\lambda_i^0+\frac{\eps}{1+\eps}\mu_i>0$. Then applying \eqref{eq:sumdiff} for $(y_i) = (y_i^0)$ and $(\lambda_i) = (\lambda_i^\eps)$,
\b*
\nabla\varphi(x_0)-h(x_0) = \sum_{i=1}^l\lambda_i^\eps c_x(x_0,y_i) = \sum_{i=1}^l\frac{1}{1+\eps}\lambda_i c_x(x_0,y_i) + \frac{\eps}{1+\eps}c_x(x_0,y).
\e*
By subtracting, we get $c_x(x_0,y)=A_{x_0}\left(\frac{1+\eps}{\eps}\sum_{i=1}^l(\lambda_i^\eps-\frac{1}{1+\eps}\lambda_i) y_i\right)=A_{x_0}(y)$. Now doing this for all $x\in\R^d$ so that $\varphi$ is differentiable in $x$, by domination of $\mu_I$ by Lebesgue, this holds for $\mu_I-$a.e. $x\in\R^d$, $\eta-$a.s. and therefore $\mu-$a.s.

\no\underline{(ii)} Now we prove the converse statement. Let $S_0\subset \{A(Y)=c_x(x_0,Y)\}$ be a closed bounded subset of $\Omega$ for some $x_0\in\R^d$, and $A\in\Aff_d$ such that $x_0\in\interior\,\conv \,S_0$, $c$ is $\Ctn^{2,0}\cap\Ctn^{1,1}$ in the neighborhood of $S_0$, and $c_{xy}(S_0)-\nabla A\subset GL_d(\R)$. First, we show that $S_0$ is finite. Indeed, we suppose to the contrary that $|S_0|=\infty$, we can find a sequence $(y_n)_{n\ge 1}\subset S_0$ with distinct elements. As $S_0$ is closed bounded, and therefore compact, we may extract a subsequence $(y_{\varphi(n)})$ converging to $y_l\in S_0$. We have
$c_x(x_0,y_{\varphi(n)}) = A(y_{\varphi(n)})$, and $c_x(x_0,y_l) = A(y_l)$. We subtract and get $c_x(x_0,y_{\varphi(n)})-c_x(x_0,y_l) - \nabla A(y_{\varphi(n)}-y_l) = 0$, and using Taylor-Young around $y_l$,
$c_{xy}(x_0,y_l)(y_{\varphi(n)}-y_l) + o(|y_{\varphi(n)}-y_l|) - \nabla A(y_{\varphi(n)}-y_l) = 0$. As $y_\varphi(n)\neq y_l$ for $n$ large enough , we may write $u_n := \frac{y_{\varphi(n)}-y_l}{|y_{\varphi(n)}-y_l|}$.
As $u_n$ stands in the unit sphere which is compact, we can extract a subsequence $(u_{\psi(n)})$, converging to a unit vector $u$. As we have
$c_{xy}(x_0,y_l)u_{\psi(n)} + o(1) - \nabla Au_{\psi(n)} = 0$,
we may pass to the limit $n\to \infty$, and get:
$$(c_{xy}(x_0,y_l)- \nabla A)u=0.$$
As $u\neq 0$, we get the contradiction: $c_{xy}(x_0,y)-\nabla A\notin GL_d(\R)$.

Now, we denote $S_0=\{y_i\}_{1\le i\le k}$ where $k:=|S_0|$. For $r>0$ small enough, the balls $\Bb\big((x_0,y_i),r\big)$ are disjoint, $c_{xy}(\cdot)-\nabla A\subset GL_d(\R)$ on these balls by continuity of the determinant, and $c$ is $\Ctn^{2,0}\cap\Ctn^{1,1}$ on these balls. Now we define appropriate dual functions. Let $M>0$ large enough so that on the balls, $(M-1)I_d-(\nabla A+\nabla A^t)-c_{xx}$ is positive semidefinite.

We set $h(X):=\nabla A(X-x_0)-A(x_0)$, and $\varphi(X):=\frac12 M|X-x_0|^2$. Now for $1\le i\le k$, $c_x(x_0,y_i)-\nabla A\cdot(y_i-x_0)=\nabla\varphi(x_0)-h(x_0)$, $(x,y)\longmapsto c_x(x,y)-\nabla A\cdot(y-x)$ is $\Ctn^1$, and its partial derivative with respect to $y$, $c_{xy}-\nabla A$ is invertible on the balls. Then by the implicit functions Theorem, we may find a mapping $T_i\in\Ctn^1(\R^d,\R^d)$ such that for $x\in\R^d$ in the neighborhood of $x_0$,
\be\label{eq:zerodiff}
c_x\big(x,T_i(x)\big)-\nabla A\cdot\big(T_i(x)-x\big)=\nabla\varphi(x)-h(x).
\ee
Its gradient at $x_0$ is given by $\nabla T_i(x_0) = \big(c_{xy}(x_0,y_i)-\nabla A\big)^{-1}\big( M I_d-(\nabla A+\nabla A^t)-c_{xx}(x_0,y_i)\big)$. This matrix is invertible, and therefore by the local inversion theorem $T_i$ is a $\Ctn^1-$diffeomorphism in the neighborhood of $x_0$. We shrink the radius $r$ of the balls so that each $T_i$ is a diffeomorphism on $B:=X\Big(\Bb\big((x_0,y_i),r\big)\Big)$ (independent of $i$). Let $B_i:=T_i(B)$, for $y\in B_i$, let $\psi(y):= c\big(T_i^{-1}(y),y\big)-\varphi\big(T_i^{-1}(y)\big)-h\big(T_i^{-1}(y)\big)\cdot\big(y-T_i^{-1}(y)\big)$. These definitions are not interfering, as we supposed that the balls $B_i$ are not overlapping.

Let $\Gamma:=\{(x,T_i(x)):x\in B,1\le i\le k\}$. By definition of $\psi$, $c=\varphi\oplus\psi+h^\otimes$ on $\Gamma$. Now let $(x,y)\in B\x B_i$, for some $i$. $(x_0,y)\in\Gamma$, for some $x_0\in B$. Let $F:= \varphi\oplus\psi+h^\otimes-c$, we prove now that $F(x,y)\ge 0$, with equality if and only if $x=x_0$ (i.e. $(x,y)\in\Gamma$). $F(x_0,y) = 0$, and $F_x(x_0,y) = 0$ by \eqref{eq:zerodiff}. However, $F_{xx}(X,Y)=MI_d-(\nabla A+\nabla A^t)-c_{xx}(X,Y)$ which is positive definite on $B\x B_i$, and therefore we get
\b*
F(x,y)=F(x,y)-F(x_0,y)&=&\int_{x_0}^x F_x(z,y)\cdot dz=\int_{x_0}^x \big(F_x(z,y)-F_x(x_0,y)\big)\cdot dz\\
&=& \int_{x_0}^x\int_{x_0}^z dw\cdot F_{xx}(w,y)\cdot dz\ge 0.
\e*
Where the last inequality follows from the fact that $F_{xx}$ is positive definite and $dw$ and $dz$ are two vectors collinear with $(x-x_0)$. It also proves that $F(x,y)=0$ if and only if $(x,y)\in\Gamma$.

Now, we define $\Ctn^1$ mappings $\lambda_i:B\longrightarrow (0,1]$ such that $\sum_{i=1}^k\lambda_i(x)T_i(x) = x$. We may do this because we assumed that $x\in\interior \,\conv \, S_0$, and therefore, by continuity, up to reducing $B$ again, $x\in\interior \,\conv\{ T_1(x),...,T_k(x)\}$ for all $x\in B$. Finally let $\mu_0\in\Pc(\R^d)$ such that $\supp\, \mu_0 = B$ with $\Ctn^\infty$ density $f$ (take for example a well chosen wavelet). Now for $1\le i\le k$, we define $\nu_0$ on $B_i$ by $\nu_0(dy) = \lambda_i\big(T_i^{-1}(y)\big)f\big(T_i^{-1}(y)\big)\left|\det \nabla T_i\big(T^{-1}(y)\big)\right|^{-1}$. Then $\P^*(dx,dy):=\mu_0(dx)\otimes\sum_{i=1}^k\lambda_i(x)\delta_{T_i(x)}(dy)$ is supported on $\Gamma$, is in $\Mc(\mu_0,\nu_0)$. As $\varphi$, and $\psi$ are continuous, and therefore bounded, and as $\mu_0$ and $\nu_0$ are compactly supported, $\P^*[c] = \mu_0[\varphi]+\nu_0[\psi]$, and therefore $\P^*$ is an optimizer for $\Sbf_{\mu_0,\nu_0}(c)$.

Now we prove that this is the only optimizer. Let $\P$ be an optimizer for $\Sbf_{\mu_0,\nu_0}(c)$. Then $\P[\Gamma] = 1$, and therefore $\P(dx,dy)=\mu_0(dx)\otimes\sum_{i=1}^k\gamma_i(x)\delta_{T_i(x)}(dy)$, for some mappings $\gamma_i$. Let $1\le i \le k$, as for $y\in B_i$, there is only one $x:=T_i^{-1}(y)\in B$ such that $(x,y)\in\Gamma$. Then we may apply the Jacobian formula: $\nu_0(dy) = \gamma_i\big(T_i^{-1}(y)\big)f\big(T_i^{-1}(y)\big)\left|\det \nabla T_i\big(T^{-1}(y)\big)\right|^{-1}$. As this density in also equal to $\nu_0(dy) = \lambda_i\big(T_i^{-1}(y)\big)f\big(T_i^{-1}(y)\big)\left|\det \nabla T_i\big(T^{-1}(y)\big)\right|^{-1}$, and as $f\big(T_i^{-1}(y)\big)\left|\det \nabla T_i\big(T^{-1}(y)\big)\right|^{-1}>0$, we deduce that $\lambda_i\big(T^{-1}(Y)\big)=\gamma_i\big(T^{-1}(Y)\big)$, $\nu_0-$a.s. and $\lambda_i = \gamma_i$, $\mu_0-$a.s. and therefore $\P=\P^*$.
\ep\\

The statement (i) of Theorem \ref{thm:structure} is well known, it is already used in \cite{henry2016explicit} (to establish Theorem 5.1), \cite{beiglboeck2016problem} (see Theorem 7.1), and \cite{ghoussoub2015structure} (for Theorem 5.5). However, the converse implication (ii) is new and we will show in the next subsections how it gives crucial information about the structure of martingale optimal transport for classical cost functions. This converse implication will serve as a counterexample generator, similar to counterexample 7.3.2 in \cite{beiglboeck2016problem}, which could have been found by an immediate application of the converse implication (ii) in Theorem \ref{thm:structure}.

Beiglb\"ock \& Juillet \cite{beiglboeck2016problem} and Henry-Labord\`ere \& Touzi \cite{henry2016explicit} solved the problem in dimension 1 for the distance cost or for costs satisfying the "Spence-Mirless condition" (i.e. $\frac{\partial^3}{\partial x\partial y^2} c>0$), in these particular cases, the support of the optimal probabilities is contained in two points in $y$ for $x$ fixed. See also Beiglb\"ock, Henry-Labord\`ere \& Touzi \cite{beiglbock2015monotone}. Some more precise results have been provided by Ghoussoub, Kim, and Lim \cite{ghoussoub2015structure}: they show that for the distance cost, the image can be contained in its own Choquet boundary, and in the case of minimization, they show that in some particular cases the image consists of $d+1$ points, which provides uniqueness. They conjecture that this remains true in general. The subsequent theorem will allow us to prove that this conjecture is wrong, and that the properties of the image can be found much more precisely.

\subsection{Algebraic geometric finiteness criterion}\label{subsect:geoalg}

\subsubsection{Completeness at infinity of multivariate polynomial families}

Algebraic geometry is the study of algebraic varieties, which are the sets of zeros of families of multivariate polynomials. When the cost $c$ is smooth, the set $\{c_x(x_0,Y) = A(Y)\}$ for $x_0\in\R^d$ and $A\in \Aff_d$, behaves locally as an algebraic variety. This statement is illustrated by Proposition \ref{prop:couplingpoly} and Theorem \ref{thm:minmapping}.

Let $k,d\in\N$ and $(P_1,...,P_k)$ be $k$ polynomials in $\R[X_1,...,X_d]$. We denote $\langle P_1,...,P_{i-1}\rangle$ the ideal generated by $(P_1,...,P_{i-1})$ in $\R[X_1,...,X_d]$ with the convention $\langle\emptyset\rangle = \{0\}$, and $P^{hom}$ denotes the sum of the terms of $P$ which have degree $\deg(P)$:
\b*
\mbox{If }P(X) = \sum_{|\alpha|\le \deg P}a_\alpha X^\alpha,&\mbox{then }&P^{hom}(X) := \sum_{|\alpha|= \deg P}a_\alpha X^\alpha.
\e*

\begin{Definition}\label{def:transversal}
Let $k,d\in\N$ and $(P_1,...,P_k)$ be $k$ multivariate polynomials in $\R[X_1,...,X_d]$. We say that the family $(P_1,...,P_k)$ is complete at infinity if
$$
QP_i^{hom}\notin \langle P_1^{hom},...,P_{i-1}^{hom}\rangle, \mbox{ for all }Q\notin \langle P_1^{hom},...,P_{i-1}^{hom}\rangle, \mbox{ for }1\leq i\leq k.\footnote{In algebraic terms this means that $P_i^{hom}$ is not a divider of zero in the quotient ring $\R[X_1,...,X_d]/\langle P_1^{hom},...,P_{i-1}^{hom}\rangle$.}
$$
\end{Definition}

\begin{Remark}\label{rmk:transversal_dim0}
This notion actually means that the intersection of the zeros of the polynomials $P_i$ in the points at infinity in the projective space has dimension $d-k-1$ (with the convention that all negative dimensions correspond to $\emptyset$), or equivalently by the correspondance from Corollary 1.4 of \cite{hartshorne2013algebraic}, that $P_1^{hom},...,P_d^{hom}$ is a regular sequence of $\R[X_1,...,X_d]$, see page 184 of \cite{hartshorne2013algebraic}. See Proposition \ref{prop:indep_infty} to understand why $P_1^{hom},...,P_d^{hom}$ may be seen as the projections of $P_1,...,P_d$ at infinity. The algebraic geometers rather say that the algebraic varieties defined by the polynomials intersect completely at infinity. The ordering of the polynomials in Definition \ref{def:transversal} does not matter. Notice that $P_1,...,P_d$ is a regular sequence if $P_1^{hom},...,P_d^{hom}$ is a regular sequence, therefore the completeness at infinity of $(P_i)_{1\le i\le k}$ implies that the intersection of the zeros of the polynomials in the points in the projective space has dimension $d-k$.
\end{Remark}

\begin{Remark}
Notice that in Definition \ref{def:transversal}, we restrict to $\R[X_1,...,X_d]$, whereas the algebraic geometry results that we will use apply with the same definition where we need to replace $\R[X_1,...,X_d]$ by $\C[X_1,...,X_d]$. However, the families $(P_i)$ that we will consider here stem from Taylor series of smooth cost functions. Therefore we only consider $(P_i)\subset \R[X_1,...,X_d]$, and we notice that in this case, Definition \ref{def:transversal} is equivalent with $\R[X_1,...,X_d]$ or with $\C[X_1,...,X_d]$, up to projecting on the real or on the imaginary part of the equations.
\end{Remark}

\begin{Example}\label{expl:monomesindep}
If $d\in \N^*$ and $k\in(\N^*)^d$ Then $(X_1^{k_1},...,X_d^{k_d})$ is complete. Indeed, let $1\leq i \leq d$, $\langle X_1^{k_1},...,X_{i-1}^{k_{i-1}}\rangle = \{X_1^{k_1}P_1 +...+ X_{i-1}^{k_{i-1}}P_{i-1},P_1,...,P_{i-1}\in\R[X_1,...,X_d]\}$. Notice that for this family of polynomials, $P\in \langle X_1^{k_1},...,X_{i-1}^{k_{i-1}}\rangle$ is equivalent to $\partial_{X^l}P(X_1 = 0,...,X_{i-1} = 0,X_i,...,X_d) = 0$ for all $l\in\N^d$ such that $l_j< k_j$ for $j<i$, and $l_j = 0$ for $j\ge i$. Let $Q\in\R[X_1,...,X_d]$ such that $QX_i^{k_i}\in\langle X_1^{k_1},...,X_{i-1}^{k_{i-1}}\rangle$, then for all such $l\in\N^d$, we have $\partial_{X^l}(QX_i^{k_i})(X_1 = 0,...,X_{i-1} = 0,X_i,...,X_d) =X_i^{k_i}\partial_{X^l}Q(X_1 = 0,...,X_{i-1} = 0,X_i,...,X_d) = 0$, and therefore $\partial_{X^l}Q(X_1 = 0,...,X_{i-1} = 0,X_i,...,X_d) = 0$, implying that $Q\in\langle X_1^{k_1},...,X_{i-1}^{k_{i-1}}\rangle$.
\end{Example}

The notion is also invariant by linear change of variables. For example, $(X^3+XY+3, Y^3-X^2+X)$ is complete at infinity because the homogeneous polynomial family $(X^3, Y^3)$ is complete at infinity by Example \ref{expl:monomesindep} above.

\begin{Example}\label{expl:gcdindep}
Let $d\in \N$ and $(P_1,P_2)$ be $2$ homogeneous polynomials in $\R[X_1,...,X_d]\setminus\R$, then $(P_1,P_2)$ is complete at infinity if and only if $\gcd(P_1,P_2)=1$. Indeed, if $\gcd(P_1,P_2)=\Pi\neq 1$, then $P_1/\Pi\notin \langle P_1\rangle$ but $P_1/\Pi P_2 = P_1P_2/\Pi\in\langle P_1\rangle$, and therefore $(P_1,P_2)$ is not complete at infinity. Conversely, if $(P_1,P_2)$ is non complete at infinity, we may find $P',Q\in \R[X_1,...,X_d]$ such that $Q\notin \langle P_1\rangle$ and $QP_2 = P_1 P'$. We assume for contradiction that $\gcd(P_1,P_2)= 1$, then $P_1$ is a divider of $Q$, and $Q\in\langle P_1\rangle$, whence the contradiction.
\end{Example}

Let $k,d\in\N$ and $(P_1,...,P_k)$ be $k$ homogeneous polynomials in $\R[X_0,X_1,...,X_d]$, we define the set of common zeros of $(P_1,...,P_k)$: $Z(P_1,...,P_k) = \{x\in\P^d:P_i(x)=0,\mbox{ for all }1\leq i \leq k\}$. An element $x\in\R^d$ is a single common root of $P_1,...,P_k$ if $x\in Z(P_1,...,P_k)$, and the vectors $\nabla P_i(x)$ are linearly independent in $\R^d$.

\begin{Remark}\label{rmk:always_transversal}
Let $k\in (\N^*)^d$. It is well known by algebraic geometers that we may find a polynomial equation system $T\in\R\big[(X_{i,j})_{1\le i\le d,j\in (\N^*)^d: |j|\le k_i}\big]$ such that for all $(P_1,...,P_d)\in \prod_{i=1}^d \R_{k_i}[X_1,...,X_d]$ with $P_i = \sum_{j\in (\N^*)^d: |j|\le k_i}a_{i,j}X_1^{j_1}...X_d^{j_d}$, we have the equivalence
$$
T\big((a_{i,j})_{1\le i\le d,j\in (\N^*)^d: |j|\le k_i}\big) \neq 0\iff (P_1,...,P_d)\mbox{ is complete at infinity}.
$$
We provide a proof of this statement in Subsection \ref{subsect:proofcardbound}. Furthermore, not all multivariate polynomials families $(P_1,...,P_d)\in \prod_{i=1}^d \R_{k_i}[X_1,...,X_d]$ are solution of $T$ as shows Example \ref{expl:monomesindep}. As a consequence $T$ is non-zero and we have that almost all (in the sense of the Lebesgue measure) homogeneous polynomial family is complete at infinity.
\end{Remark}

\subsubsection{Criteria for finite support of conditional optimal martingale transport}\label{sect:finiteness}

We start with the one dimensional case. We emphasize that the sufficient condition $(i)$ below corresponds to a local version of \cite{henry2016explicit}.

\begin{Theorem}\label{thm:bintree}
Let $d=1$ and let $S_0=\{c_x(x_0,Y)=A(Y)\}$, for some $A\in\aff(\R,\R)$, such that $x_0\in\ri\,\conv S_0$, and $c:\Omega\longmapsto \R$.

\no{\rm (i)} If $y\mapsto c_x(x_0,y)$ is strictly convex or strictly concave for some $x_0\in \R$, then $|S_0| \le 2$.

\no{\rm (ii)} If for all $y_0\in\R$, we can find $k(y_0)\geq 2$ such that $y\mapsto c_x(x_0,y_0)$ is $k(y_0)$ times differentiable in $y_0$ and $c_{xy^{k(y_0)}}(x_0,y_0)\neq 0$, then $S_0$ is discrete. If furthermore $c_x(x_0,\cdot)$ is super-linear in $y$, then $S_0$ is finite.
\end{Theorem}
\proof
\no\underline{\rm (i)} The intersection of a strictly convex or concave curve with a line is two points or one if they intersect.

\no\underline{\rm (ii)} We suppose that $S_0$ is not discrete. Then we have $(y_n)\in S_0^\N$ a sequence of distinct elements converging to $y_0\in \R$. In $y_0$, $f:y\mapsto c_x(x_0,y)$ is $k$ times differentiable for some $k\geq 2$ and $f^{(k)}(y_0)=c_{xy^k}(x_0,y_0)\neq 0$. We have $f(y_n) = A(y_n)$. Passing to the limit $y_n\to y_0$; we get $f(y_0) = A(y_0)$. Now we subtract and get $f(y_n)-f(y_0)=\nabla A(y_n-y_0)$. We finally apply Taylor-Young around $y_0$ to get
$$(f'(y_0)-\nabla A)(y_n-y_0)+\underset{i=2}{\overset{k}{\sum}}\frac{f^{(i)}(y_0)}{i!}(y_n-y_0)^i+o(|y_n-y_0|^k) = 0$$
This is impossible for $y_n$ close enough to $y_0$, as one of the terms of the expansion at least is nonzero. If furthermore $c_x(x_0,\cdot)$ is superlinear in $y$, $S_0$ is bounded, and therefore finite.
\ep

Our next result is a weaker version of Theorem \ref{thm:bintree} (i) in higher dimension.

\begin{Proposition}\label{prop:couplingpoly}
Let $x_0\in \R^d$ such that for $y\in\R^d$, $c_x(x_0,y)=\sum_{i=1}^dP_i(y)u_i$, with for $1\leq i \leq d$, $P_i\in \R[Y_1,..,Y_d]$ and $(u_i)_{1\leq i \leq d}$ a basis of $\R^d$. We suppose that the $P_i$ have degrees $deg(P_i)\geq 2$ and are complete at infinity. Then if $S_0=\{c_x(x_0,Y)=A(Y)\}$ for some $x_0\in\R^d$, and $A\in\Aff_d$, we have
$$|S_0| \leq \deg(P_1)...\deg(P_d).$$
\end{Proposition}

The proof of this proposition is reported in Subsection \ref{subsect:proofcardbound}.
\begin{Remark}\label{rmk:boundreached}
This bound is optimal as we see with the example: $P_i = (Y_i-1)(Y_i-2)...(Y_i-k_i)$, for $1\leq i\leq d$. Then $\{1,2,...,k_1\}\x...\x\{1,...,k_d\}=\{c_x(x_0,Y) = A(Y)\}$. (For $A=0$) And this set has cardinal $k_1...k_d=deg(P_1)...deg(P_d)$. But this bound is not always reached when we fix the polynomials as we can see in the example $d=1$ and $P= X^4$, we can add any affine function to it, it will never have more than $2$ real zeros even if its degree is $4$.
\end{Remark}

The following example illustrates this theorem in dimension 2.

\begin{figure}[H]
\centering

 \includegraphics[width=0.2\linewidth]{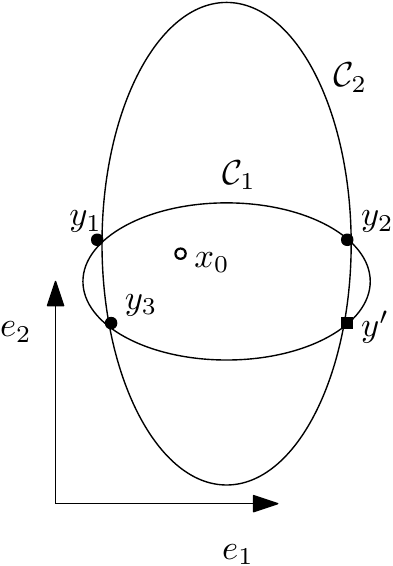}
    \caption{\label{fig:ellipses} Solution of $c_x(x_0,Y)=A(Y)$ for $c(x,y)=x_1(y_1^2+2y_2^2)+x_2(2y_1^2+y_2^2)$.}
\end{figure}

\begin{Example}\label{expl:ellipses}
Let $d=2$ and $c:(x,y)\in\R^2\x\R^2\longmapsto x_1(y_1^2+2y_2^2)+x_2(2y_1^2+y_2^2)$. Then $c_x(x,y)=(y_1^2+2y_2^2)e_1+(2y_1^2+y_2^2)e_2$ for all $(x,y)$, where $(e_1,e_2)$ is the canonical basis of $\R^2$. Let $A\in\Aff_2$, $A=A_1e_1+A_2e_2$. The equation $c_x(x_0,y)=A(y)$ can be written
\begin{equation*}
  \left\{
    \begin{aligned}
y_1^2+2y_2^2 &=& A_1(e_1)y_1+A_1(e_2)y_2+A_1(0)\\
2y_1^2+y_2^2 &=& A_2(e_1)y_1+A_2(e_2)y_2+A_2(0).
    \end{aligned}
  \right.
\end{equation*}

\no These equations are equations of ellipses $\Cc_1$ of axes ratio $\sqrt{2}$ oriented along $e_1$, and $\Cc_2$ of axes ratio $\sqrt{2}$ oriented along $e_2$. Then we see visualy on Figure \ref{fig:ellipses} that in the nondegenerate case, $\Cc_1$ and $\Cc_2$ are determined by three affine independent points $y_1,y_2,y_3\in\{c_x(x_0,Y)=A(Y)\}$, and that a fourth point $y'$ naturally appears in the intersection of the ellipses.
\end{Example}

Now we give a general result. If $k\ge 1$, we denote
\be\label{eq:notation_bracket}
c_{x_i,y^{k}}(x_0,y_0)[Y^k]&:=&\sum_{1\le j_1,...,j_k\le d}\partial^{k+1}_{x_i,y_{j_1},...,y_{j_k}}c(x_0,y_0)Y_{j_1}...Y_{j_k},
\ee
the homogeneous multivariate polynomial of degree $k$ associated to the Taylor term of the expansion of the map $c_{x_i}(x_0,\cdot)$ around $y_0$ for $1\le i\le d$.

We now provide
e the extension of Theorem \ref{thm:bintree} (ii) to higher dimension.

\begin{Theorem}\label{thm:finitenesscrit}
Let $x_0\in \R^d$ and $S_0=\{c_x(x_0,Y)=A(Y)\}$ for some $A\in\Aff_d$. Assume that for all $y_0\in \R^d$ and any $1\leq i \leq d$, $c_{x_i}(x_0,\cdot)$ is $k_i\geq 2$ times differentiable at the point $y_0$ and that $\left(c_{x_i,y^{k_i}}(x_0,y_0)[Y^{k_i}]\right)_{1\leq i \leq d}$ is a complete at infinity family of $\R[Y_1,...,Y_d]$, then $S_0$ consists of isolated points. If furthermore $c_x(x_0,\cdot)$ is super-linear in $y$, then $S_0$ is finite.
\end{Theorem}

The proof of this theorem is reported in Subsection \ref{subsect:proofcardbound}.

\subsection{Largest support of conditional optimal martingale transport plan}\label{subsect:lowerbound}

The previous section provides a bound on the cardinal of the set $S_0$ in the polynomial case, which could be converted to a local result for a sufficiently smooth function, as it behaves locally like a multivariate polynomial. However, with the converse statement (ii) of the structure Theorem \ref{thm:structure}, we may also bound this cardinality from below. Let $c$ be a $\Ctn^{1,2}$ cost function, and $x_0\in\R^d$, we denote
\b*
N_c(x_0):= \sup_{P\in\R_1[Y_{1},...,Y_{d}]^d}\big|Z^1_{\R}(H_c(x_0)+P)\big|,&\mbox{where}&H_c(x_0):= \left(c_{x_i,y^{2}}(x_0,x_0)[Y^2]\right)_{1\le i \le d}.
\e*
where we denote by $Z^1_{\R}(Q_1,...,Q_d)$ the set of real (finite) single common zeros of the multivariate polynomials $Q_1,...,Q_d\in\R[Y_1,...,Y_d]$.

\begin{Definition}
We say that $c$ is second order complete at infinity at $x_0\in \R^d$ if $c$ is differentiable at $x=x_0$ and twice differentiable at $y=x_0$, and $H_c(x_0)$ is a complete at infinity family of $\R_2[Y_1,...,Y_d]$.
\end{Definition}

\begin{Remark}
Recall that by Remark \ref{rmk:always_transversal}, this property holds for almost all cost function. We highlight here that this consideration should be taken with caution, indeed cost functions of importance which are $c:= f(|X-Y|)$ with $f$ smooth fail to be second order complete at infinity, even in the case of $c$ smooth at $(x_0,x_0)$, as the sets $\{c_x(x_0,Y) = A(Y)\}$ for $A\in\Aff_d$ may be infinite and contradict Theorem \ref{thm:finitenesscrit}, as they may contain balls, see Theorem \ref{thm:fundistance} below.
\end{Remark}

\begin{Theorem}\label{thm:minmapping}
Let $c:\Omega\longrightarrow \R$ be second order complete at infinity and $\Ctn^{2,0}\cap \Ctn^{1,2}$ in the neighborhood of $(x_0,x_0)$ for some $x_0\in\R^d$. Then, we may find $\mu_0,\nu_0\in\Pc(\R^d)$ with $\Ctn^1$ densities, and a unique $\P^*\in\Mc(\mu_0,\nu_0)$ such that
\b*
\Sbf_{\mu_0,\nu_0}(c) = \P^*[c]&\mbox{and}&|\supp\, \P^*_X| = N_c(x_0),\,\mu-\mbox{a.s.}
\e*
\end{Theorem}

The proof of this result is reported in subsection \ref{subsect:minmapping}. Theorem \ref{thm:minmapping} shows the importance of the determination of the numbers $N_c(x_0)$. We know by Remark \ref{rmk:boundreached} that for some cost $c:\Omega\longrightarrow \R$, the upper bound is reached: $N_c(x_0) = 2^d$. We conjecture that this bound is reached for all cost which is second order complete at infinity at $x_0$. An important question is whether there exists a criterion on cost functions to have the differential intersection limited to $d+1$ points, similarly to the Spence-Mirless condition in one dimension. It has been conjectured in \cite{ghoussoub2015structure} in the case of minimization for the distance cost. Theorem \ref{thm:pdistance} together with (ii) of Theorem \ref{thm:structure} proves that this conjecture is wrong. Now we prove that even for much more general second order complete at infinity cost functions, there is no hope to find such a criterion for $d$ even.

\begin{Theorem}\label{thm:mindimeven}
Let $x_0\in \R^d$, and $c$ second order complete at infinity and $\Ctn^{1,2}$ at $(x_0,x_0)$, then
\b*
d+1+\mathbf{1}_{\{d\mbox{ even}\}}&\le N_c(x_0) \le& 2^d.
\e*
\end{Theorem}

The proof of Theorem \ref{thm:mindimeven} is reported in Subsection \ref{subsect:minmapping}.

\subsection{Support of optimal plans for classical costs}\label{subsect:charactclassical}

\subsubsection{Euclidean distance based cost functions}

Theorem \ref{thm:structure} shows the importance of sets $S_0=\{c_x(x_0,Y)=A(Y)\}$ for $x_0\in\ri\,\conv\,S_0$, and $A\in\Aff_d$. We can characterize them precisely when $c:(x,y)\in\R^d\x\R^d\longmapsto f(|x-y|)$ for some $f\in\Ctn^1(\R_+,\R)$. In view of Remark \ref{rmk:d+1determine}, the following result gives the structure of $S_0$ as a function of $d+1$ known points in this set. Let $g:t> 0\longmapsto -f'(t)/t$, notice that
\b*
c_x(x,y)=g(|y-x|)(y-x),&\mbox{on}&\{X\neq Y\}.
\e*

Furthermore, $c(x,y)$ is differentiable in $x=y$ if and only if $f'(0) = 0$, in this case $c_x(x,x) = 0$. We fix $S_0:=\{ c_x(x_0,Y)=A(Y)\}$, for some $x_0\in\interior\,\conv\,S_0$, and $A\in\Aff_d$.
The next theorem gives $S_0$ as a function of $A$ and $x_0$. For $a\notin Sp(\nabla A)$, let $\mathbf{y}(a) := x_0 + (aI_d-\nabla A)^{-1}A(x_0)$. For $a\in Sp(\nabla A)$, if the limit exists, we write $|\mathbf{y}(a)|<\infty$ and denote $\mathbf{y}(a):=\lim_{t\to a}\mathbf{y}(t)$.

\begin{Theorem}\label{thm:fundistance} Let $S_0:=\{c_x(x_0,Y)=A(Y)\}$ for $x_0\in\ri\,\conv\,S_0$, and $A\in\Aff_d$. Then
\b*
S_0 = \cup_{(t,\rho)\in\Ac}S_t^\rho\cup\big\{\mathbf{y}(t):t\in \mathbf{fix}(g\circ|\mathbf{y}-x_0|)\big\},
\e*
where $S_t^\rho := \Sc_{V_t}\left(p_t,\sqrt{ \rho^2-|p_t-x_0|^2}\right)$, with $V_t := \mathbf{y}(t)+\ker(tI_d-\nabla A)$, $p_t := proj_{V_t}(x_0)$, and $\Ac:= \big\{(t,\rho) : t\in Sp(\nabla A),\,|\mathbf{y}(t)|<\infty,\,g(\rho) = t,\mbox{ and }\rho \ge|p_t-x_0|\big\}$.
\end{Theorem}

\no{\rm (i)} The elements in the spheres $S_{t_0}^\rho$ for all $\rho$ from Theorem \ref{thm:fundistance} will be said to be $2d_{t_0}$ degenerate points, where $d_{t_0} := \dim V_{t_0}$. This convention corresponds to the degree $2d_{t_0}$ of their associated root $t_0$ of the extended polynomial $\chi(t):= \det(tI_d-\nabla A)^2g^{-1}(t)^2-|Com(tI_d-\nabla A)^tA(0)|^2$). Notice that in the case $d_{t_0} = 1$, the sphere $S_{t_0}^\rho$ is a $0-$dimensional sphere, which consists in $2d_{t_0} = 2$ points.

\no{\rm (ii)} We say that $\mathbf{y}(t_0)\in S_0$ is double for $t_0\in\R$ if $\min_{t}\Big\{g\big(|\mathbf{y}(t)-x_0|\big)-t\Big\}=0$ (attained at $t_0$) where the minimum is taken in the neighborhood of $t_0$. Notice that then in the smooth case, $t_0$ is a double root of $\chi$.

\begin{Corollary}\label{corr:2d}
$S_0$ contains at least $2d$ possibly degenerate points counted with multiplicity.
\end{Corollary}

The proofs of Theorem \ref{thm:fundistance} and Corollary \ref{corr:2d} are reported in Subsection \ref{subsect:proofpdistance}.

\subsubsection{Powers of Euclidean distance cost}

In this section we provide calculations in the case where $f$ is a power function. The particular cases $p=0,2$ are trivial, for other values, we have the following theorems.

\begin{Theorem}\label{thm:pdistance}
Let $c:=|X-Y|^p$. Let $S_0:=\{ c_x(x_0,Y)=A(Y)\}$, for some $x_0\in\interior\,\conv\,S_0$, and $A\in\Aff_d$. Then if $p\le 1$, $S_0$ contains $2d$ possibly degenerate points counted with multiplicity, and if $1< p < 2-\frac25$ or $p>2+\frac23$, $S_0$ contains $2d+1$ possibly degenerate points counted with multiplicity.
\end{Theorem}

The proof of this theorem is reported in Subsection \ref{subsect:proofpdistance}.

\begin{Remark}
In both cases, for almost all choice of $y_0,...,y_d\in \R^d$ as the first elements of $S_0$, determining the Affine mapping $A$, we have $d_i = 0$ for all $i$, and $c_{xy}(x_0,S_0)-\nabla A \subset GL_d(\R^d)$. Then for $-\infty < p\leq 1$, and $p\neq 0$, $|S_0| = 2d$, and for $1< p < 2-\frac25$ or $p>2+\frac23$, $|S_0| = 2d+1$. Therefore, by (ii) of Theorem \ref{thm:structure}, we may find $\mu,\nu\in\Pc(\R^d)$ with $\Ctn^1$ densities such that the associated optimizer $\P\in\Mc(\mu,\nu)$ of the MOT problem \eqref{pb:MOT} satisfies $|\supp\,\P_X| = 2d$, $\mu-$a.s. if $p\le 1$, and $|\supp\,\P_X| = 2d+1$, $\mu-$a.s. if $p> 1$.
\end{Remark}

\begin{Remark}
Based on numerical experiments, we conjecture that the result of Theorem \ref{thm:pdistance} still holds for $2-\frac25\le p \le 2+\frac23$, and $p\neq 2$. See Section \ref{sect:numerics}.
\end{Remark}

\begin{Remark}\label{rmk:minimization}
Assumption \ref{ass:duality} implies that $c$ is subdifferentiable. Then we can deal with cost functions $c:=-|X-Y|^p$ with $0<p\le 1$ only by evacuating the problem on $\{X=Y\}$. If $0<p\le 1$, it was proved by Lim \cite{lim2014optimal} that in this case the value $\{X=Y\}$ is preferentially chosen by the problem: Theorem 4.2 in \cite{lim2014optimal} states that the mass $\mu\wedge\nu$ stays put (i.e. this common mass of $\mu$ and $\nu$ is concentrated on the diagonal $\{X=Y\}$ by the optimal coupling) and the optimization reduces to a minimization with the marginals $\mu-\mu\wedge\nu$ and $\nu-\mu\wedge\nu$. Therefore, $c$ is differentiable on all the points concerned by this other optimization, and the supports are given by $\supp\,\P_x\subset\{c_x(x,Y)=A_x(Y)\}\cup\{x\}$, for $\mu-$a.e. $x\in\R^d$. Then the supports are exactly given by the ones from the maximisation case with eventually adding the diagonal.
\end{Remark}

Notice that Remark \ref{rmk:minimization} together with (ii) of Theorem \ref{thm:structure} and Theorem \ref{thm:pdistance} prove that Conjecture 2 in \cite{ghoussoub2015structure} is wrong, and explains the counterexample found by Lim \cite{lim2016multi}, Example 2.9.

\subsubsection{One and infinity norm cost}

For $\eps\in\Ec^1:=\{-1,1\}^d$, we denote $\Qc^1_\eps:=\prod_{1\le i\le d}\eps_i(0,\infty)$ the quadrant corresponding to the sign vector $\eps$. Similarly, for $\eps\in\Ec^\infty:=\{\pm e_i\}_{1\le i\le d}$, we denote $\Qc^\infty_\eps:=\{y\in\R^d:\eps\cdot y>|y-(\eps\cdot y)\eps|_\infty\}$ the quadrant corresponding to the signed basis vector $\eps$.

\begin{Proposition}\label{prop:characone}
Let $c := |X-Y|_p$ with $p\in\{1,\infty\}$, and $S_0 := \{c_x(x_0,Y)=A(Y)\}$ for some $x_0\in\ri\,\conv\,S_0$, and $A\in\Aff_d$, with $r:=\rank\,\nabla A$. Then, we may find $2\le k\le \mathbf{1}_{p = 1}2^r+\mathbf{1}_{p = \infty}2r$, $\eps_1,...,\eps_{k}\in \Ec^p$, and $y_1,...,y_{k}\in\R^d$ such that
$$S_0= \cup_{i=1}^{k}(x_0+\Qc^p_{\eps_i})\cap(y_i+\ker\nabla A).$$

\no In particular, $S_0$ is concentrated on the boundary of its convex hull.
\end{Proposition}

This Proposition will be proved in Subsection \ref{subsect:charactoneinfty}. The case $r=d$ is of particular interest.

\begin{Remark}
Notice that the gradient of $c$ is locally constant where it exists (i.e. if $c$ is differentiable at $(x_0,y_0)$, then $c$ is differentiable at $(x,y)$ and $\nabla c(x,y) = \nabla c(x_0,y_0)$ for $(x,y)$ in the neighborhood of $(x_0,y_0)$). Then if $r=d$, $c_{xy}(x_0,S_0)-\nabla A=-\nabla A\in GL_d(\R)$, $S_0$ is finite and $|S_0| \le \mathbf{1}_{p = 1}2^d+\mathbf{1}_{p = \infty}2d$. The bound is sharp (consider for example $A:=x_0+I_d$). Therefore, by (ii) of Theorem \ref{thm:structure}, we may find $\mu,\nu\in\Pc(\R^d)$ with $\Ctn^1$ densities such that the associated optimizer $\P\in\Mc(\mu,\nu)$ of the MOT problem \eqref{pb:MOT} satisfies $|\supp\,\P_X| = \mathbf{1}_{p = 1}2^d+\mathbf{1}_{p = \infty}2d$, $\mu-$a.s.
\end{Remark}

\subsubsection{Concentration on the Choquet boundary}\label{subsect:choquet}

Recall that a set $S_0$ is included in its own Choquet boundary if $S_0\subset Ext\big(\cl{\conv}(S_0)\big)$, i.e. any point of $S_0$ is extreme in $\cl{\conv}(S_0)$. A result showed in \cite{ghoussoub2015structure} is that the image of the optimal transport is concentrated in its own Choquet boundary for distance cost. We prove that this is a consequence of (i) of the structure Theorem \ref{thm:structure}, and we generalize this observation to some other cases.

\begin{Proposition}\label{thm:choquetconv}
Let $c:\Omega\longrightarrow\R$ be a cost function, $A\in\Aff_d$, $S_0\subset\{c_x(x_0,Y)=A(Y)\}$, and $x_0\in \ri\,\conv S_0$. $S_0$ is concentrated in its own Choquet boundary in the following cases:

\no{\rm (i)} the map $y\mapsto c_x(x_0,y)\cdot u$ is strictly convex for some $u\in \R^d$;

\no{\rm (ii)} $c:(x,y)\longmapsto|x-y|_p$, with $1<p<\infty$;

\no{\rm (iii)} $c:(x,y)\longmapsto|x-y|^p$, with $-\infty<p\leq 1$;

\no{\rm (iv)} $c:(x,y)\longmapsto|x-y|^p$, with $1<p< 2-\frac25$ or $p>2+\frac23$, and $p\left(\min_{y\in S_0}|y-x_0|\right)^{p-2}$ is a double root of the polynomial $\det(\nabla A-X I_d)^2-|p|^{\frac{2}{2-p}} X^{\frac{2}{2-p}}|Com(\nabla A-X I_d)^tA(0)|^2$.

Furthermore, if $c:(x,y)\longmapsto|x-y|^p$, with $1<p< 2-\frac25$ or $p>2+\frac23$, and $S_0$ is not concentrated on its own Choquet boundary, then we may find a unique $y_0\in S_0$ such that $|y_0-x_0| = \min_{y\in S_0}|y-x_0|$, and $S_0\setminus\{y_0\}$ is concentrated on its own Choquet boundary.
\end{Proposition}

The proof of this proposition is reported in Subsection \ref{subsect:Choquetp}.

\begin{Remark}
If $p=1$ or $p=\infty$, there are counterexamples to Proposition \ref{thm:choquetconv} (ii), as $S_0$ may contain a non-trivial face of itself , see Proposition \ref{prop:characone}.
\end{Remark}

\section{Proofs of the main results}\label{sect:proof}


\subsection{Proof of the support cardinality bounds}\label{subsect:proofcardbound}

We first introduce some notions of Algebraic geometry. Recall $\P^d := \big(\C^{d+1}\big)^*/\C^*$, the $d-$dimensional projective space which complements the space with points at infinity. Recall that there is an isomorphism $\P^d\approx \C^d\cup \P^{d-1}$, where $\P^{d-1}$ are the "points at infinity". Then we may consider the points for which $x_0 = 0$ as "at infinity" because the surjection of $\P^d$ in $\C^d$ is given by $(x_0,x_1,...,x_d)\longmapsto (x_1/x_0,...,x_d/x_0)$ so that when $x_0 = 0$, we formally divide by zero and then consider that the point is sent to infinity. The isomorphism $\P^d\approx \C^d\cup \P^{d-1}$ follows from the easy decomposition:
\b*
\P^d&=& \{(x_0,...,x_d)\in \C^{d+1},x_0\neq 0\}/\C^*\cup \{(0,x_1,...,x_d),(x_1,...,x_d)\in\C^d\setminus\{0\}\}/\C^*\\
&=&\{(1,x_1/x_0,...,x_d/x_0),(x_0,...,x_d)\in \C^{d+1},x_0\neq 0\}\\
&&\cup \big\{(0,x_1,...,x_d),(x_1,...,x_d)\in\big(\C^d\big)^*\big\}/\C^*\\
&\approx& \C^d\cup\big(\C^d\big)^*/\C^*
\approx \C^d\cup \P^{d-1}.
\e*

The points in the projective space $\P^d$ in the equivalence class of $\{x_0 = 0\}$ are called points at infinity.

\begin{Definition}
The map
$$P= \sum_{n\in\N^d,|n|\leq \deg(P)}a_n X^n\longmapsto P^{proj}:=\sum_{n\in\N^d,|n|\leq \deg(P)}a_n X^nX_0^{deg(P)-|n|},$$
defines an isomorphism between $\C[X_1,...,X_d]$ and $\C^{hom}[X_0,X_1,...,X_d]$. Let $(P_1,...,P_k)$ be $k\ge 1$ polynomials in $\R[X_1,...,X_d]$, we define the set of common projective zeros of $(P_1,...,P_k)$ by $Z^{proj}(P_1,...,P_k) := Z(P_1^{proj},...,P_k^{proj})$. 
\end{Definition}

This allows us to define the zeros of a nonhomogeneous polynomial in the projective space.

We finally report the following well-known result which will be needed for the proofs of Proposition \ref{prop:couplingpoly} and Theorem \ref{thm:mindimeven}.

\begin{Theorem}[Bezout]\label{thm:Bezout}
Let $d\in\N$ and $P_1,...,P_d\in \R[X_1,...,X_d]$ be complete at infinity. Then $|Z^{proj}(P_1,...,P_d)|=deg(P_1)...deg(P_d)$, where the roots are counted with multiplicity.
\end{Theorem}
\proof
By Corollary 7.8 of Hartshorne \cite{hartshorne2013algebraic} extended to $\P^d$ and $d$ curves, we have
\be\label{eq:Bezout}
\sum_{V \in Irr\big(Z^{proj}(P_1,...,P_d)\big)} i\big(Z^{proj}(P_1),...,Z^{proj}(P_d),V\big) &=& deg(P_1)...deg(P_d),
\ee
where $i\big(Z^{proj}(P_1),...,Z^{proj}(P_d),V\big)$ is the multiplicity of the intersection of $Z^{proj}(P_1)$,..., and $Z^{proj}(P_d)$ along $V$, and $Irr\big(Z^{proj}(P_1,...,P_d)\big)$ is the collection of irreducible components of $Z^{proj}(P_1,...,P_d)$. By Remark \ref{rmk:transversal_dim0}, $Z^{proj}(P_1,...,P_d)$ has dimension $d-d = 0$ by the fact that $(P_1,...,P_d)$ is complete at infinity. Therefore, its irreducible components (in the algebraic sense) are singletons, and \eqref{eq:Bezout} proves the result.
\ep\\

Notice that we have the identity $P^{hom} = P^{proj}(X_0 = 0)$. Then $P^{hom}$ may be interpreted as the restriction to infinity of $P^{proj}$ and we deduce the following characterization of completeness at infinity that justifies the name we gave to this notion. We believe that this is a standard algebraic geometry result, but we could not find precise references. For this reason, we report the proof for completeness. For $P_1,...,P_d\in \R[X_1,...,X_d]$, we denote $Z^{\aff}(P_1,...,P_d) := \{x\in\C^d:{P}_1(x)=...={P}_d(x) = 0\}$ the set of their common affine zeros.

\begin{Proposition}\label{prop:indep_infty}
Let $P_1,...,P_d \in \R[X_1,...,X_d]$, Then the following assertions are equivalent:

\no{\rm (i)} $(P_1,...,P_d)$ is complete at infinity;

\no{\rm (ii)}  $Z^{proj}(P_1,...,P_d)$ contains no points at infinity;

\no{\rm (iii)}  $Z^\aff(P_1^{hom},...,P_d^{hom})=\{0\}$.
\end{Proposition}
\proof
We first prove $(iii)\implies (ii)$, let $x\in\P^d$ at infinity, i.e. such that $x_0 = 0$. Then by definition of the projective space, $x':=(x_1,...,x_d)\neq 0$, and by (iii) we have that $P_i^{hom}(x')\neq 0$ for some $i$. Notice that $P_i^{hom}(x') = P_i^{proj}(x)$, and therefore $P_i^{proj}(x)\neq 0$ and $x\notin Z^{proj}(P_1,...,P_d)$.

Now we prove $(i)\implies (iii)$. By definition of completeness at infinity, we have that $(P^{hom}_1,...,P^{hom}_d)$ is complete at infinity by the fact that $(P_1,...,P_d)$ is complete at infinity. By Theorem \ref{thm:Bezout}, $(P^{hom}_1,...,P^{hom}_d)$ has exactly $\deg P^{hom}_1 ... \deg P^{hom}_d$ common projective roots counted with multiplicity. However, by their homogeneity property, $(1,0,...,0)$ is a projective root of order $\deg P^{hom}_1 ... \deg P^{hom}_d$, therefore it is the only common projective root of these multivariate polynomials, in particular $0$ is their only affine common root.

Finally we prove that $(ii)\implies (i)$. In order to prove this implication, we assume to the contrary that $(i)$ does not hold. Then by Remark \ref{rmk:transversal_dim0}, we have that the dimension of this projective variety is higher than $d-(d-1) = 1$. Then we may find some $x\in Z(P^{hom}_1,...,P^{hom}_{d-1})$ which is different from $z = (1,0,...,0)$, as if $z$ was the only zero, the dimension of $Z(P^{hom}_1,...,P^{hom}_{d-1})$ would be $0$. Now we consider $x' := (0,x_1,...,x_d)\in\P^d$. As $x_0' = 0$, $x'$ is at infinity and $P^{hom}_i(x) = P^{hom}_i(x') = P^{proj}_i(x')$. Therefore, $x'\in Z(P^{proj}_1,...,P^{proj}_{d})$, contradicting (ii) by the fact that $x'$ is at infinity.
\ep\\

\no{\bf Proof of Remark \ref{rmk:always_transversal}}
Let $\Xc:= \{x\in\P^d:x_0 = 0\}$ the subset of points of $\P^d$ at infinity, and $\Yc:= H_{k_1}\x...\x H_{k_d}$, with $H_n$ the set of homogeneous polynomials of degree $n$ for $n\in\N$. The set $\Xc$ is a projective variety as the set of zeros of the polynomial $X_0$, and the set $\Yc$ is a quasi-projective variety as it is an affine space. The set $A := \{(p,P_1,...,P_d)\in \Xc\x\Yc : P_1(p) = ... = P_d(p) = 0\}$ is a set of zeros of polynomials in $\Xc\x\Yc$ (also called closed set for the Zariski topology by algebraic geometers). Notice that the set of non-complete at infinity polynomials in $\R[X_1,...,X_d]$ is exactly the projection of $A$ on $\Yc$ by Proposition \ref{prop:indep_infty}, and therefore this set is characterized by a polynomial equation system on the coefficients of the $P_i$ by Theorem 1.11 in \cite{shafarevich2013basic}, which states that the projection of closed sets for the Zariski topology in $\Xc\x\Yc$ stays closed for the Zariski topology of $\Yc$.
\ep\\

\no{\bf Proof of Proposition \ref{prop:couplingpoly}}
For $1\le i\le d$, let $A_i:= u_i\cdot A\in\aff(\R^d,\R)$. If for each $1\le i\le d$ we project this equation onto $Vect(u_i)$ along $Vect(u_j, j\neq i)$, we get:
\b*
P_i(y) = A_i(y) ,& i=1,..,d.
\e*
Thanks to the completeness at infinity of $(P_i)$, the $\widetilde{P}_i$ which are defined for $1\leq i \leq d$ by
$$\widetilde{P}_i(Z_0,Z_1,...,Z_d) := P^{proj}_i(Z_1,...,Z_d)+  \nabla A_i Z_0^{k-1}+A_i(0) Z_0^k$$
are also complete at infinity as for all $i$, we have $\widetilde{P}_i^{hom}=P_i^{hom}$. By Bezout Theorem \ref{thm:Bezout} there are $\deg(P_1)...\deg(P_d)$ common projective roots to these polynomial. These roots may be complex, infinite, or multiple, therefore the set $S_0$ which is the set of these common roots that are finite and real has its cardinal bounded by $\deg(P_1)...\deg(P_d)$.
\ep\\

\no{\bf Proof of Theorem \ref{thm:finitenesscrit}}
We suppose that $S_0$ is not discrete. Then we have $(y_n)\in S_0^\N$ a sequence of distinct elements converging to $y_0\in \R^d$. We denote $P_i(Y_1,...,Y_d) := c_{x_i,y^{k_i}}(x_0,y_0)[Y^{k_i}]$ for $1\leq i\leq d$. We know that $(P_i)_{1\leq i \leq d}$ is a complete at infinity family of $\R[Y_1,...,Y_d]$. We have $f(y_n):= c_x(x_0,y_n) = A(y_n)$. Passing to the limit $y_n\to y_0$, we get $f(y_0) = A(y_0)$. Now subtracting the terms, we get $f(y_n)-f(y_0)=\nabla A(y_n-y_0)$, and applying Taylor-Young around $y_0$, we get
\be\label{eq:taylor_zero}
(\nabla f(y_0)-\nabla A)\cdot(y_n-y_0)+\underset{i=2}{\overset{k-1}{\sum}}\frac{f^{(i)}(y_0)}{i!}[(y_n-y_0)^i]+P(y_n-y_0)+o(|y_n-y_0|^{k_i}) = 0
\ee

With $P = (P_1,...,P_d)$. By Proposition \ref{prop:couplingpoly}, the Taylor multivariate polynomial is locally nonzero around $y_0$ as it has a finite number of zeros on $\R^d$. This is in contradiction with \eqref{eq:taylor_zero} for $y_n$ close enough to $y_0$.

If furthermore $c$ is super-linear in the $y$ variable at $x_0$, $T$ is bounded, and therefore finite.
\ep

\subsection{Lower bound for a smooth cost function}\label{subsect:minmapping}

As a preparation for the proof of Theorem \ref{thm:mindimeven}, we need to prove the following lemma.

\begin{Lemma}\label{lemma:detnozero}
Let $(P_1,...,P_d)$ be a complete at infinity family in $\R_2[X_1,...,X_d]$. Then the multivariate polynomial $\det(\nabla P_1,...,\nabla P_d)$ is non-zero.
\end{Lemma}
\proof
We suppose to the contrary that $\det(\nabla P)= 0$, where we denote $P = (P_1,...,P_d)$. We claim that we may find $y_0\in\R^d$, and a map $u:\R^d\longrightarrow \Sc_1(0)$ which is $\Cc^\infty$ in the neighborhood of $y_0$ and such that $u(y)\in \ker(\nabla P(y))$ for $y$ in this neighborhood. Then we solve the differential equation $y'(t) = u(y(t))$ with initial condition $y(0) = y_0$. As a consequence of the regularity of $u$ in the neighborhood of $y_0$, by the Cauchy-Lipschitz theorem, this dynamic system has a unique solution for $t$ in a neighborhood $[-\eps,\eps]$ of $0$, where $\eps>0$. However, we notice that $P(y(t))$ is constant in $t$, indeed, $\frac{d(P(y(t)))}{dt} = \nabla P(y(t))u(y(t)) = 0$. Since $|y'(t)| = 1$, this solution is non constant, then $P-P(y_0)$ has an infinity of roots: $y([-\eps,\eps])$. However, as $P$ is non-constant, $P-P(y_0)$ is also complete at infinity, which is in contradiction with the fact that it has an infinity of zeros by the Bezout Theorem \ref{thm:Bezout}.

It remains to prove the existence of $y_0\in\R^d$, and a map $u:\R^d\longrightarrow \R^d$, $\Cc^\infty$ in the neighborhood of $y_0$, such that $u(y)\in \ker(\nabla P(y_0))$ for $y$ in this neighborhood. For all $i<d$, we consider the determinants of submatrices of $\nabla P$ which have size $i$. Let $r\ge 0$ the biggest such $i$ so that at least one of these determinants is not the zero polynomial. By the fact that $\det(\nabla P)= 0$, and that the polynomials are non-constant by completeness at infinity, we have $0<r<d-1$. We fix one of these non-zero polynomial determinants. Let $x_0\in\R^d$ such that this determinant is non-zero at $y_0$. As this determinant is continuous in $y$, it is non-zero in the neighbourhood of $y_0$. Therefore, $\nabla P$ has exacly rank $r$ in the neighbourhood of $y_0$. Now we show that this consideration allows to find a continuous map $y\longmapsto u(y)$, such that $u(y)$ is a unit vector in $\ker(\nabla P)$. Notice that $\ker(\nabla P) = \Ima(\nabla P^t)^\perp$. We consider $r$ columns of $\nabla P^t$ that are used for the non-zero determinant. We apply the Gramm-Schmidt orthogonalisation algorithm on them. We get $u_1(y),...,u_r(y)$, an orthonormal basis of $\Ima(\nabla P(y)^t)$, defined and $\Ctn^\infty$ in the neighbourhood of $y_0$. Then let $u_0\in\ker(\nabla P(y_0))$, a unit vector. The map
$$u(y):=\frac{u_0-\sum_{i=1}^r\langle u_0, u_i(y)\rangle u_i(y)}{|u_0-\sum_{i=1}^r\langle u_0, u_i(y)\rangle u_i(y)|}$$
is well defined, $\Ctn^\infty$, and in $\Ima(\nabla P(y)^t)^\perp = \ker(\nabla P(y))$ in the neigbourhood of $y_0$, and therefore satisfies the conditions of the claim.
\ep\\

\no{\bf Proof of Theorem \ref{thm:mindimeven}}

\no\underline{Step 1:} Let $P_i:=(X_1,...,X_d)c_{x_i,yy}(x_0,x_0)(X_1,...,X_d)^t$. Let $y_1,...,y_{d+1}\in\R^d$, affine independent. We may find $A\in\Aff_d$ such that $A(y_i) = P(y_i)$ for all $i$, where we denote $P:=(P_i)_{1\le i\le d}$. Now we prove that $\nabla( P(y_i')-A)$ may be made invertible at points $y_i'$ at the neighborhood of $y_i$. Recall that $A$ is a function of the $d+1$ vectors $y_i$: $A = A(y_1,...,y_{d+1})$. Then we look for an explicit expression of $\nabla A(y_1,...,y_{d+1})$ (denoted $\nabla A$ for simplicity) as a function of the $y_i$. Let $Y = Mat(y_i-y_{d+1},i=1,...,d)$, the matrix with columns $y_i-y_{d+1}$, using the equality $\nabla A y_i+A(0) = P(y_i)$, we get the identity $\nabla A Y = M$, where we denote $M:=Mat(P(y_i)-P(y_{d+1}),i=1,...,d)$. Then we get the result $\nabla A = MY^{-1}$ ($Y$ is invertible as the $y_i$ are affine independent). Then having $\nabla P(y_{d+1})-\nabla A$ invertible is equivalent to having $\nabla P(y_{d+1})Y-M$ invertible. Notice that $\nabla P(y_{d+1})Y-M = -Mat(\tilde{P}(y_i), i = 1,...,d)$, where $\tilde{P} = P-P(y_{d+1})-\nabla P(y_{d+1})\cdot(Y-y_{d+1})$, and that the multivariate polynomials $\tilde{P}_i$ are complete at infinity, as they only differ from the $P_i$ by degree one polynomials. Consider the multivariate polynomial $D:=\det(\nabla\tilde{P})$. Let $1\le i\le d$, by Lemma \ref{lemma:detnozero} we may find $y_i'$ in the neighborhood of $y_i$ such that $D(y_i')\neq 0$, and therefore $\nabla\tilde{P}(y_i')$ is invertible. Thanks to this invertibility, we may perturb the $y_i'$ to make $M':=Mat(\tilde{P}(y_i'), i = 1,...,d)$ invertible. As $Sp(M')$ is finite, for $\lambda>0$ small enough, $M'+\lambda I_d$ is invertible. For $1\le i\le d$, we may find $y_i''$ in the neighborhood of $y_i'$ so that $\tilde{P}(y_i'') = \tilde{P}(y_i')+\lambda e_i+o(\lambda)$, thanks to the invertibility of $\nabla\tilde{P}(y_i')$. Then for $\lambda$ small enough, $(P(y_i''),i=1,...,d) = M'+\lambda I_d+o(\lambda)$ is invertible.

We were able, by perturbing the $y_i$ for $i\neq d+1$ to make $\nabla( P(y_{d+1}')-A)$ invertible. By continuity, this invertibility property will still hold if we perturb again sufficiently slightly the $y_i$. Then we redo the same process, replacing $y_{d+1}'$ by another $y_i'$. We suppose that the perturbation is sufficiently small so that all the invertibilities hold in spite of the successive perturbations of the $y_i$. Finally, we found $y_1',...,y_{d+1}'$ affine independent so that $P(y_i') = A(y_i')$ and $\nabla P(y_i')-\nabla A$ is invertible for all $1\le i\le d+1$.

\no\underline{Step 2:} Then $N_c(x_0)\ge d+1$ because $y_1',...,y_{d+1}'$ are $d+1$ single real roots of $P+A = H_c(x_0)+A$, and $A\in \Aff_d$, which may be identified to $\R_1[Y_{1},...,Y_{d}]^d$. As the $P_i-A_i$ are real multivariate polynomials, all non-real zeros have to be coupled with their complex conjugate. Recall that by Theorem \ref{thm:Bezout}, there are exactly $2^d$ zeros to this system. There are no zeros at infinity by Proposition \ref{prop:indep_infty}, and there is an even number of non-real zeros by the invariance by conjugation observation. Then there must be an even number of real roots. As the $y_i'$ are simple roots by invertibility of the derivative of $P-A$ at these points, there must be an even number of real roots, counted with multiplicity. If $d$ is even, $d+1$ is odd, which proves the existence of a possibly multiple $d+2-$th zero $y_0$, distinct from the $y_i$. We assume, up to renumbering, that $y_0',...,y_d'$ are affine independent, and we perturb again $y_0',...,y_d'$ to make $y_0$ a single zero. We need to check that $y_{d+1}'$ is still a single zero of $P-A$. Indeed, the map $(y_1',...,y_{d+1}')\longmapsto A$ if locally a diffeomorphism around $(y_1,...,y_{d+1})$, then by the implicit functions Theorem, we may write $y_{d+1}' = F(y_1',...,y_{d}',A) = F\big(y_1',...,y_{d}',A(y_0',...,y_d')\big)$, where $F$ is a local smooth function. Then $y_{d+1}'$ remains a single zero if the perturbation of $y_0,...,y_d$ is small enough. The result is proved, if $d$ is even we may find $d+2$ single zeros to $P-A$.

The reverse inequality is a simple application of Proposition \ref{prop:couplingpoly}.
\ep\\

As a preparation for the proof of Theorem \ref{thm:minmapping}, we introduce the two following lemmas:

\begin{Lemma}\label{lemma:interiorZ1}
Let $Q_1,...,Q_d$, $d$ complete at infinity multivariate polynomials of degree $2$ and $x\in\R^d$. Then, for all $P_1,...,P_d$ multivariate polynomials of degree $1$, we may find $\widetilde{P}_1,...,\widetilde{P}_d$, multivariate polynomials of degree $1$ such that $|Z^1_\R(Q_1+\widetilde{P}_1,...,Q_d+\widetilde{P}_d)| \ge |Z^1_\R(Q_1+P_1,...,Q_d+P_d)|$ and $x\in \interior \,\conv\,Z^1_\R(Q_1+\widetilde{P}_1,...,Q_d+\widetilde{P}_d)$.
\end{Lemma}
\proof
Let $P_1,...,P_d$ multivariate polynomials of degree $1$. We claim that we may find $R_1,...,R_d$ of degree $1$ so that $Z^1_\R(Q_1+R_1,...,Q_d+R_d)$ has full dimension and contains $Z^1_\R(Q_1+P_1,...,Q_d+P_d)$. Then we may find $x'\in\interior\,\conv\,Z^1_\R(Q_1+R_1,...,Q_d+R_d)$, and by the fact that all $Q_i$ have degree $2$, we may find $\widetilde{P}_1,...,\widetilde{P}_d$ of degree $1$ such that $(Q+P)(X+x'-x) = Q+\widetilde{P}$. Finally, as the change of variables $X+x'-x$ does not change the number of roots of $Q+P$ nor their multiplicity, and by the fact that $x\in\interior\,\conv\,Z^1_\R(Q_1+\widetilde{P}_1,...,Q_d+\widetilde{P}_d)$ by translation, $\widetilde{P}_1,...\widetilde{P}_d$ solves the problem.

Now we prove the claim. We prove by induction that we may add dimensions to $Z^1_\R(Q_1+R_1,...,Q_d+R_d)$ by changing the $R_i$. First by Theorem \ref{thm:mindimeven}, we may assume that $Z^1_\R(Q_1+R_1,...,Q_d+R_d)$ is non-empty. Up to making a distance-preserving linear change of variables, we may assume that $Z^1_\R(Q_1,...,Q_d)\subset \{X_d = 0\}$ and that $0\in Z^1_\R(Q_1,...,Q_d)$. We look for $D\in\R[X_1,...,X_d]^d$ in the form $D = X_d v$ for some $v\in\R^d$, so that $Q+D$ leaves $Z^1_\R(Q_1,...,Q_d)$ unchanged. In order to include some $y\in \{X_d\neq 0\}$, we set $D:= -Q(y)/y_dX_d$. The constraint that we have now is to fix $y$ is that $\nabla (Q+D)(y')\in GL_d(\R)$ for all $y'\in Z^1_\R(Q_1,...,Q_d)$ and for $y' = y$. Notice that all these constraints have the form $\det\big(\nabla(y_dQ-Q(y)X_d)(y')\big)\neq 0$ if $y'\neq y$, and $\det\big(\nabla(y_dQ-Q(y)X_d)(y)\big)\neq 0$ for the case $y' = y$, therefore in all the cases this is a polynomial equation in $y$. We claim that each of these equations on $y$ have a solution. Then as there is a finite number of such equations, the set of solutions is a dense open set, in particular it is non-empty and we may find $y\in \R^d$ so that $\{y\}\cup Z^1_\R(Q)\subset Z^1_\R(Q+D)$ and $\dim Z^1_\R(Q+D)>\dim Z^1_\R(Q)$. By induction, we may reach full dimension for $\dim Z^1_\R(Q+D)$, and the problem is solved.

Finally, we prove the claim that the solution set to $\det\big(\nabla(y_dQ-Q(y)X_d)(y')\big)\neq 0$ is non-empty.

\no\underline{Case 1:} $y'\in Z^1_\R(Q)$. Then, up to applying a translation change of variables, we may assume that $y' = 0$. Then by the fact that $Q$ has degree $2$, the equation that we would like to satisfy is
$$\det\Big(y_d\nabla Q(0)-\big(\nabla Q(0)y+\frac12 D^2Q(0)[y^2]\big)e_d^t\Big)\neq 0.$$
We make it more tractable by making operations on the columns:
\b*
&&\det\Big(y_d\nabla Q(0)-\big(\nabla Q(0)y+\frac12 D^2Q(0)[y^2]\big)e_d^t\Big)\\
&=& \det\left(y_d\sum_{i=1}^d\nabla Q_i(0)e_i^t-\left(\sum_{i=1}^d\nabla Q_i(0)y_i+\frac12 D^2Q(0)[y^2]\right)e_d^t\right)\\
&=& \det\left(y_d\sum_{i=1}^{d-1}\nabla Q_i(0)e_i^t-\frac12 D^2Q(0)[y^2]e_d^t\right),
\e*
where we have subtracted the $i^{th}$ column multiplied by $y_i/y_d$ to the $d^{th}$ column for all $1\le i \le d-1$. Now we prove that this multivariate polynomial is non-zero. We assume for contradiction that it is zero. Then for all $y\in\{X_d\neq 0\}$, $D^2Q(0)[y^2]\in H:= Vect(\nabla Q_i(0),1\le i\le d-1)$, which is $d-1-$dimensional by the fact that $\nabla Q\in GL_d(\R)$ by simplicity of the root $0$. By continuousness, we have in fact that $D^2Q(0)[y^2]\in H$ for all $y\in\R^d$. Therefore, for all $y_1,y_2\in\R^d$, we have the equality $D^2Q(0) [y_1,y_2] = \frac12\left(D^2Q(0) [(y_1+y_2)^2]-D^2Q(0) [y_1,y_1] -D^2Q(0) [y_2,y_2]  \right)\in H$. Then we may find $u\in\R^d$ non-zero such that $\sum_{i=1}^du_i D^2Q_i(0) = 0$. Then $(Q^{hom}_1,...,Q^{hom}_d)$ is $d-1-$dimensional and $Z^{proj}(Q_1,...,Q_d)$ is at least $1-$dimensional, then it intersects the variety of points at infinity, which is a contradiction by Proposition \ref{prop:indep_infty} together with the fact that $(Q_1,...,Q_d)$ is a complete at infinity family.

\no\underline{Case 2:} $y' = y$. Then the equation that we would like to satisfy is
$$\det\Big(y_d\nabla Q(y)-\big(\nabla Q(0)y+\frac12 D^2Q(0)[y^2]\big)e_d^t\Big)\neq 0,$$
which may be expanded thanks to the fact that $Q$ has degree $2$:
$$\det\Big(y_d\left(\nabla Q(0)+D^2Q(0)y\right)-\big(\nabla Q(0)y+\frac12 D^2Q(0)[y^2]\big)e_d^t\Big)\neq 0.$$
Similar than in the previous case, by the same operations on the columns we get:
\b*
&&\det\Big(y_d\left(\nabla Q(0)+D^2Q(0)y\right)-\big(\nabla Q(0)y+\frac12 D^2Q(0)[y^2]\big)e_d^t\Big)\\
&=& \det\left(y_d\sum_{i=1}^d\left(\nabla Q_i(0)+D^2Q_i(0)y\right)e_i^t-\left(\sum_{i=1}^d\nabla Q_i(0)y_i+\frac12 D^2Q_i(0)y y_i\right)e_d^t\right)\\
&=& \det\left(y_d\sum_{i=1}^{d-1}\left(\nabla Q_i(0)+D^2Q_i(0)y\right)e_i^t+\frac12 D^2Q(0)[y^2]e_d^t\right),
\e*
Now we assume for contradiction that this polynomial in $y$ is zero. Then for all $y\in\{X_d\neq 0\}$ small enough so that $\nabla Q(0)+D^2 Q(0)y\in GL_d(\R)$, $D^2Q(0)[y^2]\in H_y:= Vect(\nabla Q_i(0)+D^2 Q_i(0)y,1\le i\le d-1)$. Notice that up to multiplying $y$ by $\lambda>0$, we have that $\lambda^2 D^2Q(0)[y^2]\in H_{\lambda y}$, and therefore $D^2Q(0)[y^2]\in H_{\lambda y}$. By passing to the limit $\lambda \longrightarrow 0$, we have $D^2Q(0)[y^2]\in H_{0}$ thanks to the fact that $\nabla Q\in GL_d(\R)$. Therefore we obtain a contradiction similar to case $1$.
\ep\\

\begin{Lemma}\label{lemma:globinv}
Let $M>0$, we may find $R(M)$ such that for all $F:\R^d\longmapsto \R^d$ and $x_0\in\R^d$ such that on $B_{M^{-1}}(x_0)$, $F$ is $\Ctn^2$ and we have that $\nabla F$ and $D^2F$ is bounded by $M$, and $\det\nabla F\ge M^{-1}$, we have that $F$ is a $\Ctn^1-$diffeomorphism on $B_{R(M)}(x_0)$.
\end{Lemma}
\proof
The determinant is a polynomial application, therefore it is Lipschitz when restricted to the compact of matrices bounded by $M$. Let $L(M)$ be its Lipschitz constant. Then on the neighbourhood $B_{R_0(M)}(x_0)$, we have that $\det \nabla F$ is bigger than $\frac12 M^{-1}$, with $R_0(M) = \min\left(M^{-1},\frac{1}{2L(M)M}\right)$. We claim that $F$ is injective on $B_{R_1(M)}(x_0)$ with $R_1(M) := \min\left(M^{-1},\frac{1}{4M^2C(M)}\right)$, where $C(M)$ is a bound for the comatrices of matrices dominated by $M$. Then by the global inversion theorem, $F$ is a $\Ctn^1-$diffeomorphism on $B_{R(M)}(x_0)$ with $R(M) = \min\big(R_0(M),R_1(M)\big)$.

Now we prove the claim that $F$ is injective on $B_{R_1(M)}(x_0)$. Let $x,y\in B_{R_1(M)}(x_0)$,
\b*
F(y)-F(x) &=& \int_0^1\nabla F(tx+(1-t)y)(y-x)dt\\
&=& \nabla F(x)(y-x)+\int_0^1\int_0^t D^2F(sx+(1-s)y)[(y-x)^2]dsdt\\
&=& \nabla F(x) \left(y-x+\nabla F(x)^{-1}\int_0^1 (1-s)D^2F(sx+(1-s)y)[(y-x)^2]ds\right).
\e*
Then we assume that $F(y) = F(x)$. Then
\be\label{eq:squaredomination}
|y-x| &=& \left|\nabla F(x)^{-1}\int_0^1 (1-s)D^2F(sx+(1-s)y)[(y-x)^2]ds\right|\nonumber\\
&\le& \|\nabla F(x)^{-1}\|\frac{M}{2}|y-x|^2\nonumber\\
&\le& C(M) M^2|y-x|^2,
\ee
where the last estimate comes from the comatrix formula \eqref{eq:comatrix}. Then by the fact that $R_1(M)\le \frac{1}{4M^2C(M)}$, we have $|y-x|< \frac{1}{M^2C(M)}$, and therefore $x=y$ by \eqref{eq:squaredomination}. The injectivity is proved.
\ep\\

\no{\bf Proof of Theorem \ref{thm:minmapping}}
By Taylor expansion of $c_x$ in $y$ in the neighborhood of $x_0$, we get for $h\in\R^d$ and $\eps>0$ small enough that
$$c_x(x_0,x_0+\eps h)=c_x(x_0,x_0) +c_{xy}(x_0,x_0)\eps h + Q(\eps h) +  \eps^2R_\eps(h),$$
where, recalling the notation \eqref{eq:notation_bracket}, $Q_i(Y) := \frac12 c_{x_iyy}(x_0,x_0)[Y^2]$ and the remainder
$$R_\eps(h) = \int_0^1(1-t)\big(c_{xyy}(x_0,x_0 + \eps th)-c_{xyy}(x_0,x_0)\big)[h^2]dt.$$
Notice that $\nabla R_\eps(h) = 3\int_0^1(1-t)\big(c_{xyy}(x_0,x_0 + \eps th)-c_{xyy}(x_0,x_0)\big)[h]dt$.
By Proposition \ref{prop:couplingpoly}, we see that $N_c(x_0)$ is finite by second order completeness at infinity of $c$ at $(x_0,x_0)$. We consider from the definition of $N_c(x_0)$ an affine map $A\in\Aff_d$ such that the $d-$tuple of multivariate polynomials of degree one $A(X_1,...,X_d)$ satisfies
\b*
\left|Z_\R^1\big(Q_i+A(X_1,...,X_d)_i:1\le i\le d\big)\right|&=n&:=N_c(x_0).
\e*
By Theorem By Lemma \ref{lemma:interiorZ1}, let $P=(P_1,...,P_d)$, $d$ multivariate polynomials of degree $1$ such that $|Z^1_\R(Q_1+\widetilde{P}_1,...,Q_d+\widetilde{P}_d)| \ge n$ and $0\in \interior \,\conv\,Z^1_\R(Q_1+\widetilde{P}_1,...,Q_d+\widetilde{P}_d)$.

Let $A_\eps(y) := -\eps^2P(0)-\eps \nabla P(y-x_0)+ c_x(x_0,x_0) + c_{xy}(x_0,x_0)(y-x_0)$, we have that $c_x(x_0,x_0+\eps h) = A_\eps(x_0+\eps h)$ and $c_{xy}(x_0,x_0+\eps h)\in GL_d(\R)$ if and only if $Q(h)+P(h)+ R_\eps(h) = 0$ and $(\nabla Q+\nabla P)(h)+\nabla R_\eps(h)\in GL_d(\R)$.

Now let $h_1,...,h_n\in\R^d$ the $n$ elements of $Z^1_\R(Q_1+\widetilde{P}_1,...,Q_d+\widetilde{P}_d)$. By continuousness of $c_{xyy}$ in the neighborhood of $(x_0, x_0)$, up to restricting to a compact neighborhood, $c_{xyy}$ is uniformly continuous on this neighborhood. For $\eps>0$ small enough, each $x_0+\eps h_i$ in in the interior of this neighborhood. Therefore, by uniform continuousness $R_\eps$, and $\nabla R_\eps$ converges uniformly to $0$ when $\eps \longrightarrow 0$. Let $1\le i\le n$, we have $(Q+P+R_\eps)(h_i) = R_\eps(h_i)$, and $\nabla (Q+P)(h_i)\in GL_d(\R)$ by the fact that $h_i$ is a single root of $Q+P$, and therefore $\nabla (Q+P+R_\eps)(h_i)\in GL_d(\R)$ for $\eps$ small enough. Therefore we may apply Lemma \ref{lemma:globinv} around $h_i$: $Q+P+R_\eps$ is a diffeomorphism in a neighborhood of $h_i$ depending only on the lower bounds of $\det \nabla (Q+P+R_\eps)(h_i)$ and of the bounds for $\nabla (Q+P+R_\eps)$ and $D^2 (Q+P+R_\eps)$, which may then work for all $\eps$ small enough. Then for $\eps$ small enough, we may find $h_i^\eps$ in this neighborhood of $h_i$ such that $(Q+P+R_\eps)(h_i^\eps) = 0$. Furthermore, by the fact that $\nabla (Q+P+R_\eps)(h_i)\longrightarrow \nabla (Q+P)(h_i)$ when $\eps\longrightarrow 0$, $|h_i^\eps-h_i|\le 2\|\nabla (Q+P)^{-1}(h_i)\||R_\eps(h_i)|$, and therefore $h_i^\eps \longrightarrow h_i$ when $\eps\longrightarrow 0$. Then for $\eps$ small enough, the $h_i^\eps$ are distinct, $0\in\ri\,\conv(y_i^\eps,1\le i\le n)$, $Q(h_i^\eps)+P(h_i^\eps)+ R_\eps(h_i^\eps) = 0$, and $(\nabla Q+\nabla P)(h_i^\eps)+\nabla R_\eps(h_i^\eps)\in GL_d(\R)$.

Now the theorem is just an application of (ii) of Theorem \ref{thm:structure} to $S_0 := \{x_0+\eps h_i^\eps, i=1,...,n\}$.
\ep

\subsection{Characterization for the p-distance}\label{subsect:charactoneinfty}

Fot $p \geq 1$ and $x\in \R^d$, we have $c(\cdot,y)$ differentiable on $(\R^d)^*$ with
$$c_x(x,y)=\frac{1}{|x-y|^{p-1}_p}\underset{i=1}{\overset{d}{\sum}}|x_i-y_i|^{p-1}\frac{x_i-y_i}{|x_i-y_i|}e_i$$
For $p = 1$ and $p=\infty$, it takes a simpler form.

If $p=1$, $c(\cdot,y)$ is differentiable on $\prod_{i = 1}^d (\R\setminus\{y_i\})$ and $c_x(x,y) = \underset{i=1}{\overset{d}{\sum}}\frac{x_i-y_i}{|x_i-y_i|}e_i$.

If $p=\infty$, $c(\cdot,y)$ is differentiable on $\{x'\in \R^d,|x_i'-y_i|>|x_j'-y_j|,j\neq i, \text{ for some }1\leq i \leq d\}$, let $i:={\rm argmax}_{1\leq j \leq d}(|x_j-y_j|)$, we have $c_x(x,y) = \frac{x_i-y_i}{|x_i-y_i|}e_i$.


\no{\bf Proof of Proposition \ref{prop:characone}}
We start with the case $p=1$. We suppose without loss of generality that $x_0=0$. Recall that $c(\cdot,y)$ is differentiable on $(\R^*)^d$ and $c_x(0,y) = \sum_{i=1}^d\frac{y_i}{|y_i|}e_i$. Then the equation that we get is $A(y)=\sum_{i=1}^d sg(y_i)e_i$. Let $E:=\{\sum_{i=1}^d sg(y_i)e_i:y\in S_0\}\subset \eps\in\{-1,1\}^d$. We have $E\subset \Ima A$, which is an affine space of dimension $r$. Then there are $r$ coordinates $i_1,...,i_r$ that can be chosen arbitrarily in $\Ima A$, and the other coordinates are affine functions of the previous one. We denote $I:=(i_1,...,i_r)$ and $\Ib:=(1,...,d)\setminus I$. Thus, ${\rm card} (\Ima A\cap\{-1,1\}^d) \le {\rm card}(\{-1,1\}^I)=2^r$. As $0\in\ri S_0$, $r\ge 1$. Now for all $\eps\in E$, let $y_\eps\in S_0$ such that $c_x(0,y_\eps)=\eps$. Then if $y:=y_\eps + y_0\in \Qc^1_\eps$ with $y_0\in\ker \nabla A$, we have $A(y)=c_x(0,y)$, and therefore $y\in S_0$, proving the first part of the result.

Now we prove that $S_0\subset \partial \conv\,S_0$. Let us suppose to the contrary that $y\in\ri\,\conv\,S_0\cap S_0$. Let $y_1,...,y_n\in S_0$ such that $y=\sum_{i=1}^n\lambda_iy_i$, convex combination. Then $c_x(0,y) = \sum_{i=1}^n\lambda_i c_x(0,y_i)$. As $|c_x(0,y)| = \sum_{i=1}^n\lambda_i |c_x(0,y_i)| = \sqrt{d}$, we are in a case of equality in Cauchy-Schwartz inequality. $\eps:=c_x(0,y),c_x(0,y_1),...,c_x(0,y_n)$ are all non-negative multiples of the same unit vector, and therefore all equal as they have the same norm. Then $y,y_1,...,y_n\in\Qc^1_\eps$, and $y,y_1,...,y_n\in y_\eps+\ker \,\nabla A$. As we may apply the same to any $y'\in y_\eps+\ker \,\nabla A$, these vectors cannot be written as convex combinations of elements of $S_0$ that do not belong to $y_\eps+\ker \,\nabla A$. Therefore, $(y_\eps+\ker \,\nabla A)\cap S_0 =(y_\eps+\ker \,\nabla A)\cap \Qc^1_\eps$ is a face of $\conv\,S_0$. As we assumed that $y\in\ri\,\conv\,S_0$, we have $(y_\eps+\ker \,\nabla A)\cap \Qc^1_\eps=\ri\,\conv\,S_0$, by the fact that $\ri\,\conv\,S_0$ and $(y_\eps+\ker \,\nabla A)\cap \Qc^1_\eps$ are faces of $\conv\,S_0$ (which constitute a partition of $\conv\,S_0$, see Hiriart-Urruty-Lemar\'echal \cite{hiriart2013convex}) both containing $y$. This is impossible as $0\in\ri\,\conv\,S_0$ and $0\notin \Qc^1_\eps$. Whence the required contradiction.

The proof of the case $p=\infty$ is similar to the proof of Proposition \ref{prop:characone}, replacing by ${\rm card}(\{-1,1\}(e_i)_{1\le i\le d}) = 2d$ instead of $2^d$, and by  $|c_x(0,y)| = 1$ instead of $\sqrt{d}$.
\ep

\subsection{Characterization for the Euclidean p-distance cost}\label{subsect:proofpdistance}

By the fact that $\interior\,\conv\,S_0$ contains $x_0$, we may find $y_1,...,y_{d+1}\in S_0$ that are affine independent. Then we may find unique barycenter coefficients $(\lambda_i)_i$ such that $x_0=\sum_{i=1}^{d+1}\lambda_iy_i$. For some $y_1,...,y_{d+1}\in S_0$. For all $a\in\R$, we define
\be\label{eq:def_y_and_G}
\mathbf{y}'(a):=G(a)\sum_{i=1}^{d+1}\frac{\lambda_i}{a-a_i}y_i,&\mbox{with}&G(a) = \left(\sum_{i=1}^{d+1}\frac{\lambda_i}{a-a_i}\right)^{-1},\text{ and }a_i:=g(|y_i-x_0|),
\ee
where $\{b_1,...,b_r\} := \{a_1,...,a_{d+1}\}$ with $r\leq d+1$ and $b_1<...<b_r$, and $d_i := \big|\{j:a_j = b_i\}\big|-1$, the multiplicity of each $b_i$ for all $i$.

\begin{Proposition}\label{prop:linkAy}
We have $\mathbf{y}'(a) = \mathbf{y}(a)$ for all $a\notin Sp(\nabla A)$. In particular the map $\mathbf{y}'$ is independent of the choice of $y_1,...,y_{d+1}\in S_0$. Furthermore, $G(a) = \frac{(a-a_1)...(a-a_{d+1})}{\det(a I_d-\nabla A)}=\frac{(a-b_1)...(a-b_{r})}{(a-\gamma_1)...(a-\gamma_{r-1})}$ where $\gamma_1<...<\gamma_{r-1}$ are eigenvalues of $\nabla A$. Finally if we have $x_0\in \interior\,\conv(y_1,...,y_{d+1})$, then we have $b_1<\gamma_1<b_2<...<\gamma_{r-1}<b_r$.
\end{Proposition}
\proof
We suppose that $x_0 = 0$ for simplicity. Let $a\notin Sp(\nabla A)$, $\mathbf{y}(a)$ is the unique vector such that
\be\label{eq:basicequation}
\left(a I_d-\nabla A\right)\mathbf{y}(a) = A(0)
\ee
We now find the barycentric coordinates of $\mathbf{y}(a)$. For any $i$, $A(y_i) = a_iy_i$ with $a_i := g(|y_i|)$. As $(y_i)_i$ is a barycentric basis, we may find unique $(\lambda_i(a))_i\subset \R$ such that $\mathbf{y}(a) = \sum_i\lambda_i(a)y_i$, and $1=\sum_i\lambda_i(a)$. Then we apply $A$ and get $A(\mathbf{y}(a)) = \sum_i\lambda_i(a)A(y_i)$, so that $a \mathbf{y}(a) = \sum_i\lambda_i(a)a_i y_i$. Subtracting the previous equality on $\mathbf{y}(a)$, we get $0=\sum_i\lambda_i(a)(a-a_i) y_i$. As $(y_i)_i$ is a barycentric basis, it is a family or rank $d$. Then, by the fact that $\sum_{i=1}^{d+1}\lambda_i y_i = 0$, we have $(\lambda_i)_{1\le i \le d+1}$ and $(\lambda_i(a)(a-a_i))_{1\le i \le d+1}$ are in the same $1-$dimensional kernel of the matrix $(y_{1},...,y_{d+1})$. Then we may find $G(a)$ such that $\lambda_i(a)(a-a_i)=G(a)\lambda_i$. Now we assume that $a$ is not part of the $a_i$, then we have $\lambda_i(a)=G(a)\frac{\lambda_i}{a-a_i}$, and $G(a)=\left(\sum_{i=1}^{d+1}\frac{\lambda_i}{a-a_i}\right)^{-1}$. Finally
\be\label{eq:exprbary}
\mathbf{y}(a)=\mathbf{y}'(a)=G(a)\sum_{i=1}^{d+1}\frac{\lambda_i}{a-a_i}y_i&\mbox{ with }&G(a)=\left(\sum_{i=1}^{d+1}\frac{\lambda_i}{a-a_i}\right)^{-1}.
\ee

Now we prove that $G(a) = \frac{(a-a_1)...(a-a_{d+1})}{\det(a I_d-\nabla A)}$. We first assume that $a_1<...<a_{d+1}$ and that $x_0\in\interior\,\conv(y_1,...,y_{d+1})$ (i.e. $\lambda_1,...,\lambda_{d+1}>0$). Then $G(a)^{-1}$ has $d+1$ single poles $a_1,...,a_{d+1}$, such that $\lim_{a\uparrow a_i}G(a)^{-1} = +\infty$, and $\lim_{a\downarrow a_i}G(a)^{-1} = -\infty$ for all $i$. Therefore, $G(\gamma_i)^{-1} =0$ for some $a_{i}<\gamma_i<a_{i+1}$ for all $i\le d$. Then $\gamma_i$ is a pole of $G$, and $|\mathbf{y}'(a)|$ goes to infinity when $a\to\gamma_i$, as the coefficient in the affine basis $(y_i)_i$ go to $\pm\infty$. Therefore, $\gamma_i$ is an eigenvalue of $\nabla A$, as there are $d$ such eigenvalues, we have obtained all of them. Finally, by the fact that the rational fraction $f$ has degree $1$, as the set of its roots is restricted to the $d+1$ numbers $a_i$. Furthermore the $\gamma_i$ are $d$ poles, and $a^{-1}G(a)\longrightarrow (\sum_{i=1}^{d+1}\lambda_i)^{-1}=1$, when $a\to\infty$, we deduce the rational fraction $G(X) = \frac{(X-a_1)...(X-a_{d+1})}{(X-\gamma_1)...(X-\gamma_{d})}=\frac{(X-a_1)...(X-a_{d+1})}{\det(XI_d-\nabla A)}$.

Now if we chose other affine independent $(y_i)_{1\le i\le d+1}$ (this time not necessary with $x_0\in\conv(y_i,1\le i\le d+1)$), let the associated barycenter coordinates $\lambda_1,...,\lambda_{d+1}\in\R^*$, we suppose that the $(a_i)_i$ are still distinct, the poles of $\mathbf{y}'(a)$ are still the $d$ distinct eigenvalues of $\nabla A$ that are determined by the $\gamma_i$ such that $\lim_{a\to\gamma_i}|\mathbf{y}(a)|$, independent of the choice of $(y_i)_i$ because $\mathbf{y}'(a) = \left(a I_d-\nabla A\right)^{-1} A(0)$ is independent of this choice. However, the numerator of the fraction can be determined in the same way than it is determined in the previous case.

Now we want to generalize this result to $\lambda_1,...,\lambda_{d+1}\in \R$, and any $(a_i)_i$. If we stay in the open set in which $(y_i)_i$ is an affine basis of $\R^d$, the mapping $(y_i,a_i)_i\longmapsto A$ is continuous, and so is the mapping $(y_i)_i\longmapsto (\lambda_i)_i$. Therefore, as $(y_i,a_i,\lambda_i)_i\longmapsto \sum_{i=0}^d\frac{\lambda_i}{X-a_i}$ is continuous as well, the identity remains true for all $a_i,y_i$ such that $(y_i)_i$ is an affine basis and $\lambda_i\ge 0$.

Let us now focus on the multiple $a_i$s. We consider $1\leq i\leq r$ such that $d_i>0$. By passing to the limit $n\to\infty$ with some distinct $a_i^n$ converging to $a_i$ for all $1\le i \le d$, $d_i$ eigen values of $\nabla A$ at least will be trapped between the $a_i$s, as $a_i^n<\gamma_{i+1}^n<a_{i+1}^n<...<\gamma_{i+k}^n<a_{i+k}^n$ becomes at the limit $a_i=\gamma_{i+1}=a_{i+1}=...=\gamma_{i+k}=a_{i+k}$. Now we prove that no other eigenvalue is equal to $a_i$. Indeed, rewriting \eqref{eq:exprbary} that equation become
\be\label{eq:exprbaryb}
\mathbf{y}(a)=\mathbf{y}'(a)=G(a)\sum_{i=1}^r\frac{\lambda_i'}{a-b_i}y_i&\mbox{ with }&G(a)=\left(\sum_{i=1}^r\frac{\lambda_i'}{a-b_i}\right)^{-1}.
\ee
with $\lambda_i':=\sum_{a_j = b_i}\lambda_j$. And $G(a) =\frac{(X-b_1)^{d_1+1}...(X-b_{r})^{d_r+1}}{\det(XI_d-\nabla A)}$. By a similar reasoning than when the $(a_i)_i$ are distinct, we may find $b_1<\gamma_1<b_2<...<\gamma_{r-1}<b_r$, eigenvalues of $\nabla A$. Then, as $\deg \det(XI_d-\nabla A) =d$, and $(X-b_1)^{d_1}...(X-b_{r})^{d_r}$ is a divider to $\det(XI_d-\nabla A)$, we have $\det(XI_d-\nabla A)=(X-\gamma_1)...(X-\gamma_{r-1})(X-b_1)^{d_1}...(X-b_{r})^{d_r}$.
\ep 

\begin{Remark}
Notice that in Proposition \ref{prop:linkAy}, the eigenvalues of $\nabla A$ are given by the $\gamma_i$, and by each $b_i$ such that $d_i>0$, which has multiplicity $d_i$, in particular, these coefficients (up to their numbering) do not depend on the choice of $y_1,...,y_{d+1}$.
\end{Remark}

\no {\bf Proof of Theorem \ref{thm:fundistance}}
We suppose again that $x_0 = 0$ for simplicity. We know that if $y\in S_0$, $c_x(0,y) = g(|y|)y = A(y)$. We denote $a := g(|y|)$ and get,
\be\label{eq:basicequation2}
\left(a I_d-\nabla A\right)y = A(0)
\ee

Let $a\in\mathbf{fix}(g\circ|\mathbf{y}-x_0|)$, then $(aI_d-\nabla A)\mathbf{y}(a)=A(0)$, and $A\big(\mathbf{y}(a)\big)=a\mathbf{y}(a)=g\big(|\mathbf{y}(a)|\big)\mathbf{y}(a)=c_x\big(0,\mathbf{y}(a)\big)$, and therefore $\mathbf{y}(a)\in S_0$. Conversely, if $y\in S_0$ and $a:=g(|y|)$ is not an eigenvalue of $\nabla A$, $y=(aI_d-\nabla A)^{-1}A(0) = \mathbf{y}(a)$, and finally $g(|\mathbf{y}(a)|)=a$, hence $a\in\mathbf{fix}(g\circ|\mathbf{y}-x_0|)$.

Now let $t\in Sp(\nabla A)$ such that $|\mathbf{y}(t)|<\infty$. Let $y\in S_t^{\rho}$, we have $(tI_d-\nabla A)y = (tI_d-\nabla A)(y-\mathbf{y}(t))+A(0)=A(0)$, by passing to the limit $a\longrightarrow t$ in the equation $(aI_d-\nabla A)\mathbf{y}(a) = A(0)$. Finally, as $|y|^2 = \sqrt{ \rho^2-|p_t|^2}^2+|p_t|^2 = \rho^2$ by Pythagoras theorem, $A(y) = c_x(0,y)$, and therefore $y\in S_0$. Conversely, if $y\in S_0$ with $g(|y|) = t$, then we have $y-\mathbf{y}(t)\in \ker(tI_d-\nabla A)$, and $|y-p_t| = \sqrt{ \rho^2-|p_t|^2}$ by Pythagoras theorem: by definition $y\in S_t^\rho$.
\ep\\

\no {\bf Proof of Corollary \ref{corr:2d}}
We use the notations from Proposition \ref{prop:linkAy} and assume that $x_0\in\interior\,\conv(y_1,...,y_{d+1})$. By Theorem \ref{thm:fundistance}, $S_0$ contains $2\sum_{i=1}^rd_i$ degenerate points. Furthermore, for all $1\le i\le r-1$, $\lim_{t\to\gamma_i}|\mathbf{y}(t)-x_0|=\infty$, therefore, as $b_{i+1}$ is a root of $g\big(|\mathbf{y}(t)-x_0|\big)-t$ between $\gamma_i$ and $\gamma_{i+1}$, there is another root $b_i'$, possibly multiple equal to $b_i$, by continuity of $g$. Finally we have $2\sum_{i=1}^rd_i+r+(r-2) = 2d$ elements in $S_0$ at least, with possible degeneracy.
\ep\\

\no{\bf Proof of Theorem \ref{thm:pdistance}} We assume again that $x_0 = 0$ for simplicity. We suppose again that $x_0 = 0$ for simplicity. By identity \eqref{eq:comatrix}, if we multiply \eqref{eq:basicequation} by the comatrix, we get $\det(\lambda I_d-\nabla A) y = Com(\lambda I_d-\nabla A)^tA(0)$. Now taking the square norm, we get: $\det(\lambda I_d-\nabla A)^2|p|^{\frac{2}{p-2}}\lambda^{\frac{2}{p-2}} - |Com(\lambda I_d-\nabla A)^tA(0)|^2 = 0$. The polynomial with real exponents $\chi:=\det(H-X I_d)^2 - |p|^{\frac{2}{2-p}}\lambda^{\frac{2}{2-p}}|Com(X I_d-\nabla A)^tA(0)|^2$ is continuous in $(y_i)_i$, then similar to the proof of Theorem \ref{thm:fundistance}, we can pass to the limit from sequences of $y_i^n$ converging to $y_i$ for all $i$ such that for all $n\ge 1$, the vectors $y_i^n$ have distinct norms. It follows that $b_i$ is a $d_i$-eigenvalue of $\nabla A$, and a $(2d_i-1)$-root of $\chi$. By Theorem \ref{thm:fundistance}, we have
$$S_i = \Sc_{V_i}\left(p_i,\sqrt{ b_i^2-|p_i|^2}\right)\subset \{c_x(0,Y)=A(Y)\}.$$
With the radius $\sqrt{ b_i^2-|p_i|^2}>0$ as there are more than one elements in the sphere. We have a single sphere as the function $g$ is monotonic, and therefore injective.

Now we prove that if $-\infty < p \leq 1$, then the polynomial with real exponents
\be\label{eq:def_chi1}
\chi(X) := \det(XI_d-\nabla A)^2 -|p|^{\frac{2}{2-p}} X^{\frac{2}{2-p}}|Com(XI_d-\nabla A)^tA(0)|^2
\ee

has exactly $2d$ positive roots, counted with multiplicity. By Corollary \ref{corr:2d}, it has at least $2d$ roots, counted with multiplicity. Now we prove that there are at most $2d$ roots.

By Theorem \ref{thm:fundistance}, the roots of $\det(XI_d-\nabla A)$ all have the same sign (same than $p$). Consequently, the coefficients of $\det(XI_d-\nabla A)$ are alternated or all have the same sign. The same happens for $\det(XI_d-\nabla A)^2$. Now we use the Descartes rule\footnote{The Descartes rule states that for a polynomial with possibly non integer real coefficients, the number of positive roots is dominated by the number of alternations of signs of its coefficients ordered by their associated exponents, see \cite{jameson2006counting}.} for polynomials with non integer exponents in order to dominated the number of roots of $\chi$. Recall that $\chi = \det(XI_d-\nabla A)^2-|p|^{\frac{2}{2-p}} X^{\frac{2}{2-p}}|Com(XI_d-\nabla A)^tA(0)|^2$. We saw that the coefficients from the part $\det(XI_d-\nabla A)^2$ are alternated or all of the same sign. The exponent sequences from $\det(XI_d-\nabla A)^2$, and from $|p|^{\frac{2}{2-p}} X^{\frac{2}{2-p}}|Com(XI_d-\nabla A)^tA(0)|^2$ have both integer differences between two exponents from the same sequence. Then the exponents of $|p|^{\frac{2}{2-p}}X^{\frac{2}{2-p}}|Com(XI_d-\nabla A)^tA(0)|^2$ are located between the ones of $\det(XI_d-\nabla A)^2$ in the exponent sequence of $\chi$, i.e. the sequence of $\chi$ consists in one exponent from $\det(XI_d-\nabla A)^2$, then one exponent from $|p|^{\frac{2}{2-p}}X^{\frac{2}{2-p}}|Com(XI_d-\nabla A)^tA(0)|^2$, and so on. By the fact that $deg(\det(XI_d-\nabla A)^2) = 2d$ and $deg(|Com(XI_d-\nabla A)^tA(0)|^2) = 2d-2$, and $0< \frac{2}{2-p}\leq 2$. Then $\chi(X)$ has at most $2d$ alternations in its coefficients, and therefore it has at most $2d$ positive roots according to the Descartes rule.

Now, assume that $1< p < 2-\frac25$ or $p>2+\frac23$, then
\be\label{eq:def_chi2}
\chi(X) := \det(XI_d-\nabla A)^2 - |p|^{\frac{2}{2-p}}X^{\frac{2}{2-p}}|Com(XI_d-\nabla A)^tA(0)|^2
\ee
has exactly $2d+1$ positive roots counted with multiplicity.

Let us first prove that the polynomial has less than $2d+1$ roots. Similar to above, the coefficients of $\det(XI_d-\nabla A)$ are alternated. And the same happens for $\det(XI_d-\nabla A)^2$. Using the Descartes rule for polynomials with non integer coefficients, by the fact that the coefficients of $|p|^{\frac{2}{2-p}}X^{\frac{2}{2-p}}|Com(XI_d-\nabla A)^tA(0)|^2$ are located between the ones of $\det(XI_d-\nabla A)^2$, except strictly less than $3$, and as $deg(\det(XI_d-\nabla A)^2) = 2d$, it follows that $deg(|Com(XI_d-\nabla A)^tA(0)|^2) = 2d-2$ and $-3 < \frac{2}{2-p}< 5$. Then $\chi(X)$ has at most $2d+2$ alternations in its coefficients by the same reasoning than the case $p\le 1$. Furthermore, the sign of the coefficients in front of the extreme monomials are opposed (because $\chi$ is a difference of positive polynomials) then the maximum number of positive roots is odd, and therefore it has at most $2d+1$ positive roots according to Descartes rule.

By Corollary \ref{corr:2d}, we have $2d$ elements in $S_0$, more precisely, which range between $b_1$ and $b_r$. Furthermore, between $0$ and $b_1$ we can find some $a\in D$:

\no\underline{Case 1:} We assume that $p>2$. Then $\chi(X)\to -\infty$ when $X\to 0$ as we have that $- |p|^{\frac{2}{2-p}}X^{\frac{2}{2-p}}|Com(XI_d-\nabla A)^tA(0)|^2$ becomes dominant.

\no\underline{Case 2:} We assume that $p<2$. Then $\chi(X)\to -\infty$ when $X\to +\infty$ as  we have that $- |p|^{\frac{2}{2-p}}X^{\frac{2}{2-p}}|{}^tCom(XI_d-\nabla A)A(0)|^2$ becomes dominant.

Therefore there is one more real root, on the side where the polynomial goes to $-\infty$ as there is already one.
Finally $\chi$ has $2d+1$ roots at least and less than $2d+1$ roots, it follows that it has exactly $2d+1$ roots. We proved the second part of the theorem.
\ep

\subsection{Concentration on the Choquet boundary for the p-distance}\label{subsect:Choquetp}

\no{\bf Proof of Proposition \ref{thm:choquetconv}}
\no\underline{(i)} Let $y_0,y_1,...,y_k\in S_0$ such that $y_0 = \underset{i=1}{\overset{k}{\sum}}\lambda_i y_i$, convex combination. Then as $c_x(x_0,y_i)\cdot u = {}^tu A(y_i-x_0)$, we have $
\sum_{i=1}^k\lambda_i c_x(x_0,y_i)\cdot u = u^tA(y_0-x_0)= c_x(x_0,y_0)\cdot u$. As $y\mapsto c_x(x_0,y)\cdot u$ is strictly convex, this imposes that $\lambda_i = 1$ and $y_i = y_0$ for some $i$. Finally, $y_0$ is extreme in $S_0$, $S_0$ is concentrated in its own Choquet boundary.

\no\underline{(ii)} We know that for any $y\in S_0$ we have $c_x(x_0,y)=A(y)$. As the situation is invariant in $x_0$, we will assume $x_0 = 0$ for notations simplicity. We consider $1<q<+\infty$ such that $\frac{1}{p}+\frac{1}{q} = 1$. For any $y\in (\R^d)^*$,

\b*
|c_x(0,y)|_q &=& \left|\frac{1}{|y|^{p-1}_p}\underset{i=1}{\overset{d}{\sum}}|y_i|^{p-1}\frac{y_i}{|y_i|}e_i\right|_q
=\frac{1}{|y|^{p-1}_p}\left(\underset{i=1}{\overset{d}{\sum}}|y_i|^{(p-1)q}\right)^{\frac{1}{q}}
=\frac{1}{|y|^{\frac{p}{q}}_p}|y|_p^{\frac{p}{q}}
= 1,
\e*
as we know that $y\neq 0$ because $c$ is superdifferentiable. Then for any $y\in S_0$, we have $|Hy+v|_q=1$. We now assume that $y_0=\overset{k}{\underset{i=1}{\sum}}\lambda_iy_i$ is a strict convex combination with $(y_i)_{0\leq i\leq k}\in S_0^{k+1}$.
\b*
1&=&|A(y_0)|_q
=\left|\sum_{i=1}^k\lambda_i\big(A(y_i)\big)\right|_q
\leq  \sum_{i=1}^k\lambda_i|A(y_i)|_q
=\sum_{i=1}^k\lambda_i
= 1
\e*
We are in a case of equality for the triangular inequality for the norm $|\cdot|_q$. We know then that all the $\lambda_i A(y_i)$ and $A(y_0)$ are positively multiples. As we know that all their $q$-norm is $\lambda_i\neq 0$ and $1$, therefore $A(y_0) = ... = A(y_k)$ and $\frac{1}{|y_0|^{p-1}_p}\underset{i=1}{\overset{d}{\sum}}|(y_0)_i|^{p-1}\frac{(y_0)_i}{|(y_0)_i|}e_i = ... = \frac{1}{|y_k|^{p-1}_p}\underset{i=1}{\overset{d}{\sum}}|(y_k)_i|^{p-1}\frac{(y_k)_i}{|(y_k)_i|}e_i$. Notice that for $y\in\R^d$, we have $\frac{1}{|y|^{p-1}_p}\underset{i=1}{\overset{d}{\sum}}|y_i|^{p-1}\frac{y_i}{|y_i|}e_i = f(y/|y|_p)$, where $f:y\longmapsto \underset{i=1}{\overset{d}{\sum}}|y_i|^{p-1}\frac{y_i}{|y_i|}e_i$ is bijective $\R^d\longrightarrow \R^d$ for $p>1$. Then we have $\frac{y_0}{|y_0|_p} = ... = \frac{y_k}{|y_k|_p}$. It means that they all belong to the same semi straight line originated in $0$. As we supposed that $y_0$ is not extreme, $0$ can be included in the convex combination as we must have $1\leq i\leq k$ such that $|y_k|>|y_0|$. Then increasing the corresponding $\lambda_i$ while decreasing all the others, $0$ can be included. As $0\in \ri\,\conv\,S_0$, we can then put any element of $S_0$ in the convex combination and $S_0\subset \{0\}+\frac{y_0}{|y_0|}\R_+$. As $0\in \ri\,\conv\,S_0$, then $S_0= \{0\}$ and $y_0 = 0$, which is the required contradiction because we supposed that $y_0$ is not extreme in $S_0$.

\no\underline{(iii)} We use the notations from Theorem \ref{thm:pdistance}. We suppose again without loss of generality that $x_0=0$. Let $d:= \dim S_0$, for any $y_1,...,y_{d+1}\in S_0$ with full dimension $d$, we may find unique barycentric coordinates $(\lambda_i)_{1\le i\le d+1}$ such that $\overset{d}{\underset{i=0}{\sum}}\lambda_i y_i = 0$. Let $y\in S_0$ such that $p|y|^{p-2} = g(|y|)\notin Sp(\nabla A)$. By Proposition \ref{prop:linkAy}, $y$ can be expressed as
\b*
y = G(X)\sum^{d+1}_{i=1}\frac{\lambda_i}{X-a_i } y_i&\mbox{with}&G(X) = \left(\sum^{d+1}_{i=1}\frac{\lambda_i}{X-a_i}\right)^{-1}.
\e*
with $X= p|y|^{p-2} >0$. To have $y\in \conv(S_0)$ we then need to have all the $\frac{\lambda_i}{X-a_i}$ of the same sign. As we supposed that the $(a_i)_i$ is an increasing sequence, there must be a $0\leq i_0 \leq d-1$ such that $\lambda_i <0$ if $i\leq i_0$ and $\lambda \geq 0$ if $i\geq i_0+1$ (or $\lambda_i >0$ if $i\leq i_0$ and $\lambda \leq 0$ if $i\geq i_0+1$ but we will only treat the first case as this one can be dealt with similarly). Then the idea consists in proving that $\chi$ defined by \eqref{eq:def_chi1} has no zero in $]a_{i_0},a_{i_0+1}[$.

First let us prove that $G$ has no pole on $]a_{i_0},a_{i_0+1}[$. $G^{-1}$ can hit $0$ at most $d$ times (It is a polynomial of degree $d$ divided by another polynomial). It hits $0$ in any $]a_{i},a_{i+1}[$ for $i\neq i_0$, as the limits on the bounds are $+\infty$ and $-\infty$. This provides $d-1$ zeros. If there where a zero in $]a_{i_0},a_{i_0+1}[$, it would be double, as the infinity limits at $a_{i_0}^+$ and $a_{i_0+1}^-$ have the same sign. Which would be a contradiction.

Finally, as the poles of $G$ are the eigenvalues of $\nabla A$ and do not depend on the choice of $y_1,...,y_{d+1}$, we know that there are exactly two roots of $\chi$ between two poles. As $a_{i_0}$ and $a_{i_0+1}$ are two zeros surrounded by two consecutive poles, there are not other zeros between these two poles. $\chi$ has no zero on $]a_{i_0},a_{i_0+1}[$.

If $X = a_{i_0}$ or $X=a_{i_0+1}$, then it is a zero of $a_{i_0}-X$, and all the elements in the convex combination have same size than $y$. By the fact that we are in the case of equality in the Cauchy-Schwartz inequality, this proves that the combination only contains one element. Hence, $y\in S_0$ has to be extreme in $S_0$.

Now if $y$ corresponds to an eigenvalue of $\nabla A$, let $b:=g(|y|)$. We suppose that $y=\sum_{i=1}^{d+1}\mu_i y_i$, convex combination with $y_1,...,y_{d+1}\in S_0$, affine basis. Recall that all $\mathbf{y}(a)$ for $a\notin Sp(\nabla A)$ can be written $\mathbf{y}(a) = G(a)\sum_{i=1}^{d+1}\frac{\lambda_i}{a-a_i} y_i=G(a)\sum_{i=1}^{r}\frac{\lambda_i'}{a-b_i} y_i'$ where $\lambda_i' = \sum_{a_j = b_i}\lambda_j$, and $y_i' = \sum_{a_j = b_i}\frac{\lambda_j}{\lambda_i'}y_j$. Let $i_0$ such that $b_{i_0} = b$, let $\{y_1',...,y_{d_{i_0}}'\}:=\{y'\in\{y_1,...,y_{d+1}\}:g(|y'|)=b_{i_0}\}$. $y\in\aff(y_1',...,y_{d_{i_0}}')$, therefore $\mu_i=0$ if $a_i\neq b$. As $S_i$ is a sphere, it is concentrated on its own Choquet boundary, and therefore the convex combination $y=\sum_{i=1}^{d+1}\mu_i y_i$ is trivial, $y=y_i$ for some $i$ and $\mu_i=1$.

\no\underline{(iv)} In the first case, if $p|y_0|^{p-2}$ is a double root of $\chi$ defined by \eqref{eq:def_chi2}, then if $p < 2-\frac25$ or $p>2+\frac23$, $\chi$ has $2d+1$ roots and at most $2d$ distinct roots set around the poles of $G$ in the same way than in the case $p\leq 1$ in the proof of \rm{(iii)}.

The same happens when we remove the smallest element $y_0$ of $S_0$. Similarly $S_0\setminus\{y_0\}$ is concentrated on its own Choquet boundary.

Now we prove that $S_0$ is not concentrated on its own Choquet boundary. If $p|y_0|^{p-2}$ is a single root of $\chi$, we select $y_1',...,y_{d+1}'\in S_0$ such that $0$ is in their convex hull. By Proposition \ref{prop:linkAy}, if $y\in S_0$ and $X:=p|y|^{p-2}$, then
\be\label{eq:barybary}
y = G(X)\sum^{d+1}_{i=1}\frac{\lambda_i}{X-a_i } y_i
&\mbox{with}&
G(X) = \left(\sum^{d+1}_{i=1}\frac{\lambda_i}{X-a_i }\right)^{-1}.
\ee

\no\underline{Case 1:} We assume that $y_1' = y_0$. Then we apply \eqref{eq:barybary} to $X:=p|y|^{p-2}$ the second smallest zero of $\chi$ which is strictly smaller than the first pole by Theorem \ref{thm:pdistance} (which also means that $G(X)\geq 0$): $y := G(X)\sum^d_{i=0}\frac{\lambda_i}{a_i - X} y_i\in S_0$, or written otherwise:
$$\frac{\lambda_0 G(X)}{X-a_0}y_0 = G(X)\sum^{d+1}_{i=2}\frac{\lambda_i}{X-a_i } y_i -y$$
$G$ has its first zero at $a_0$ which is smaller than its first pole which is between $a_1$ and $a_2$ strictly, so that $G(X)>0$. This gives the result, rewriting the barycenter equation, we get:
$$y_0 = \sum^{d+1}_{i=2}\frac{\lambda_i(X-a_0)}{\lambda_0(X-a_i )} y_i +\frac{G(X)}{\lambda_0}y$$
Therefore, $y_0\in \conv(S_0\setminus\{y_0\})$.

\no\underline{Case 2:} Now we assume that $y_0' \neq y_0$. We write the barycenter equation for $X=p|y_0|^{p-2}$, we get:
\b*
y_0 = \sum^d_{i=0}\frac{\lambda_iG(X)}{X-a_i } y_i'
&\mbox{with}&
G(X) = \left(\sum^d_{i=0}\frac{\lambda_i}{X-a_i }\right)^{-1}.
\e*
Then for any $i$, $\frac{\lambda_iG(X)}{X-a_i}>0$ as all the $\frac{\lambda_i}{X-a_i }$ have the same sign. Therefore $y_0\in \conv(S_0\setminus\{y_0\})$.
\ep

\section{Numerical experiment}\label{sect:numerics}

\begin{figure}[H]
\centering
 \includegraphics[width=0.45\linewidth]{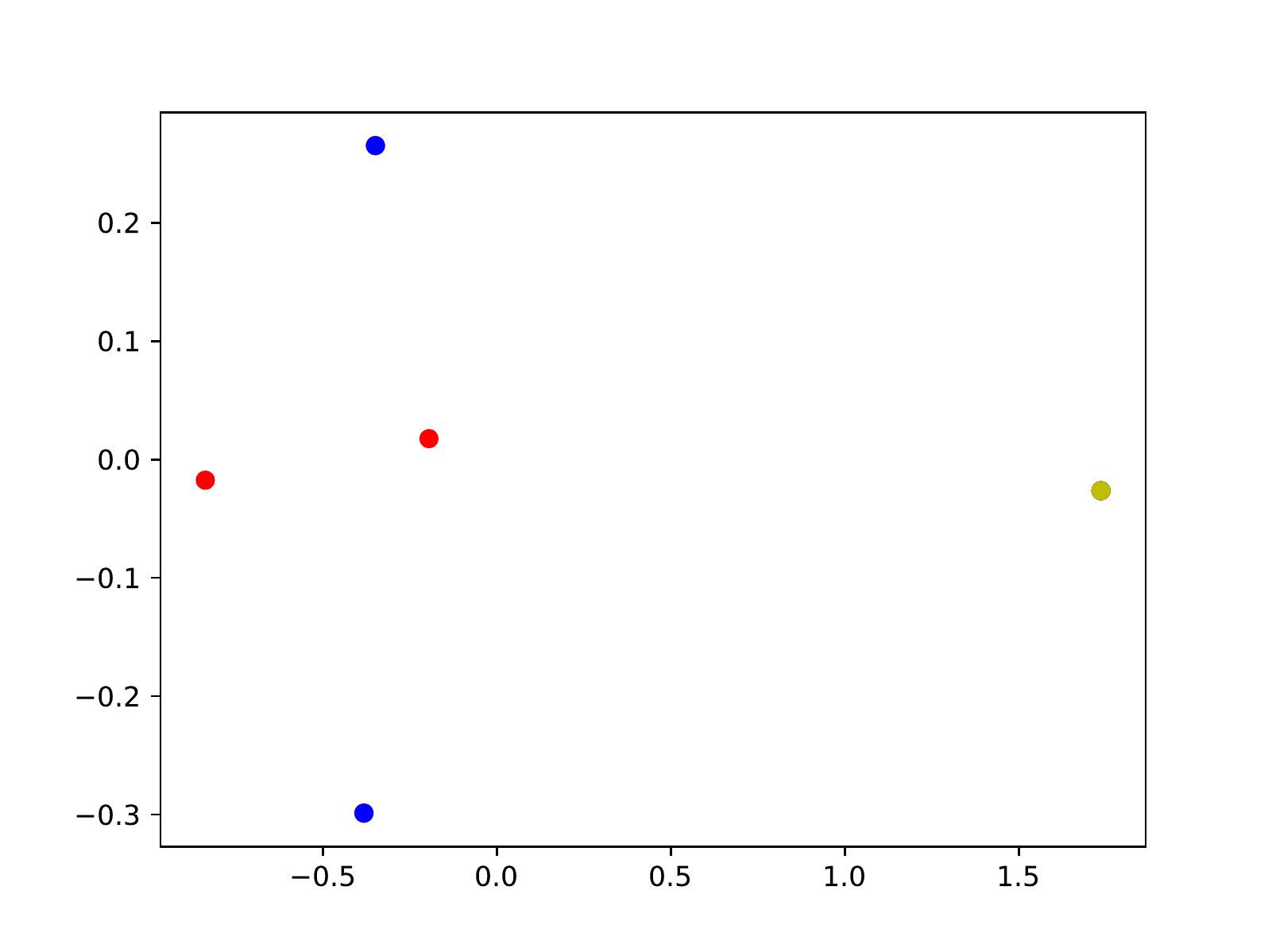}\includegraphics[width=0.45\linewidth]{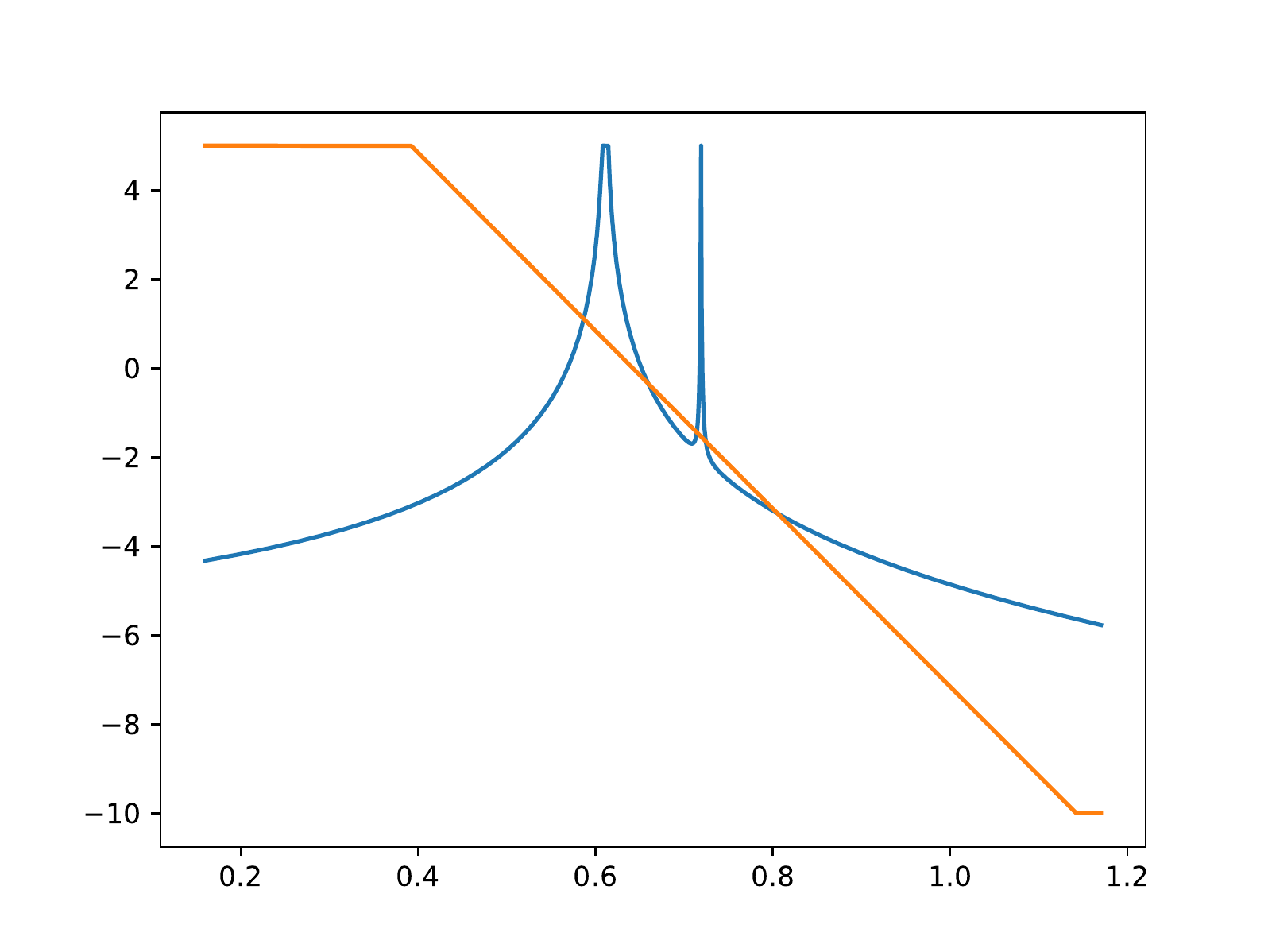}
    \caption{\label{fig:p1.9d2}$S_0$ for $d=2$ and $p=1.9$.}
\end{figure}

\begin{figure}[H]
\centering
 \includegraphics[width=0.45\linewidth]{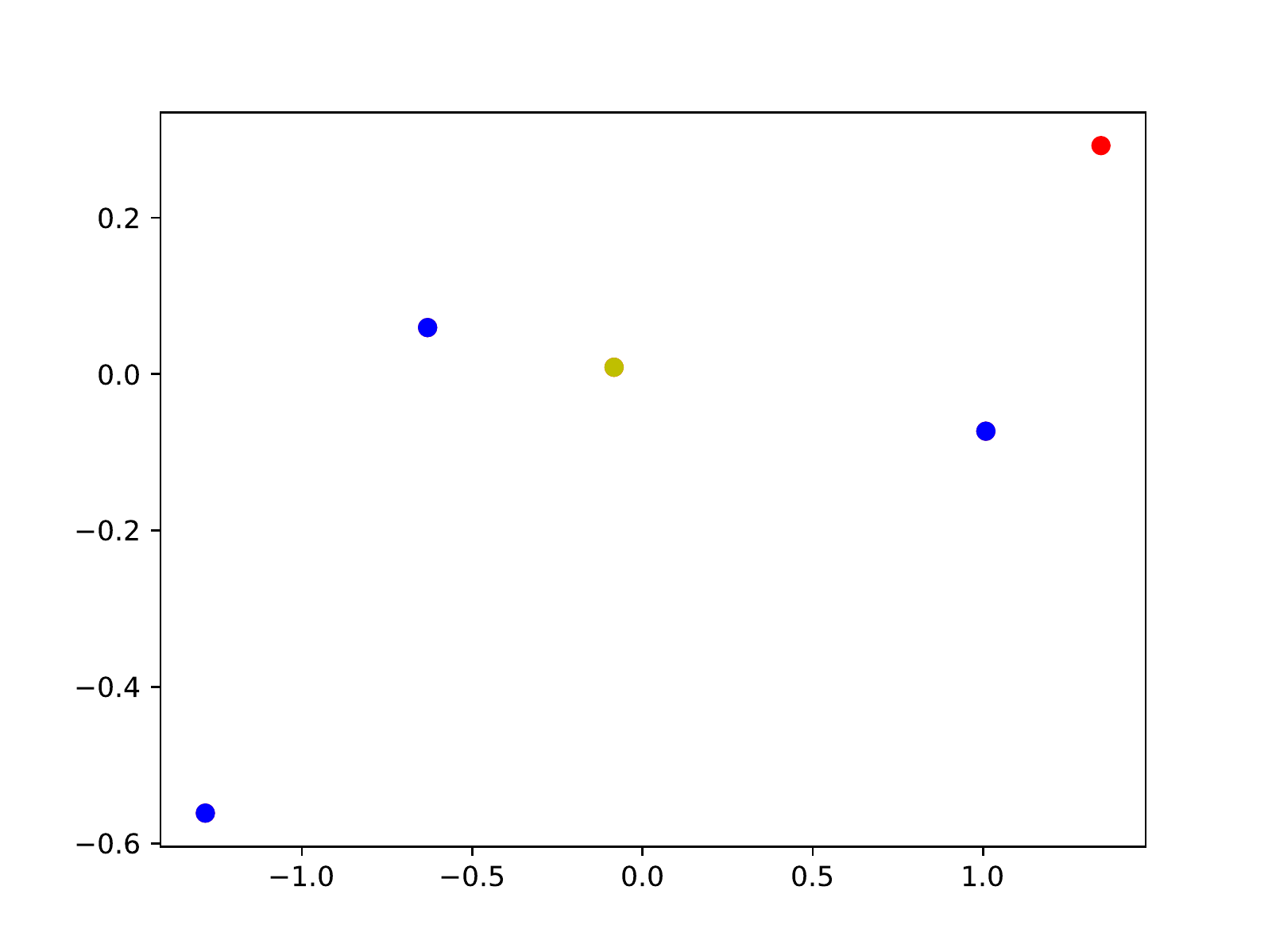}\includegraphics[width=0.45\linewidth]{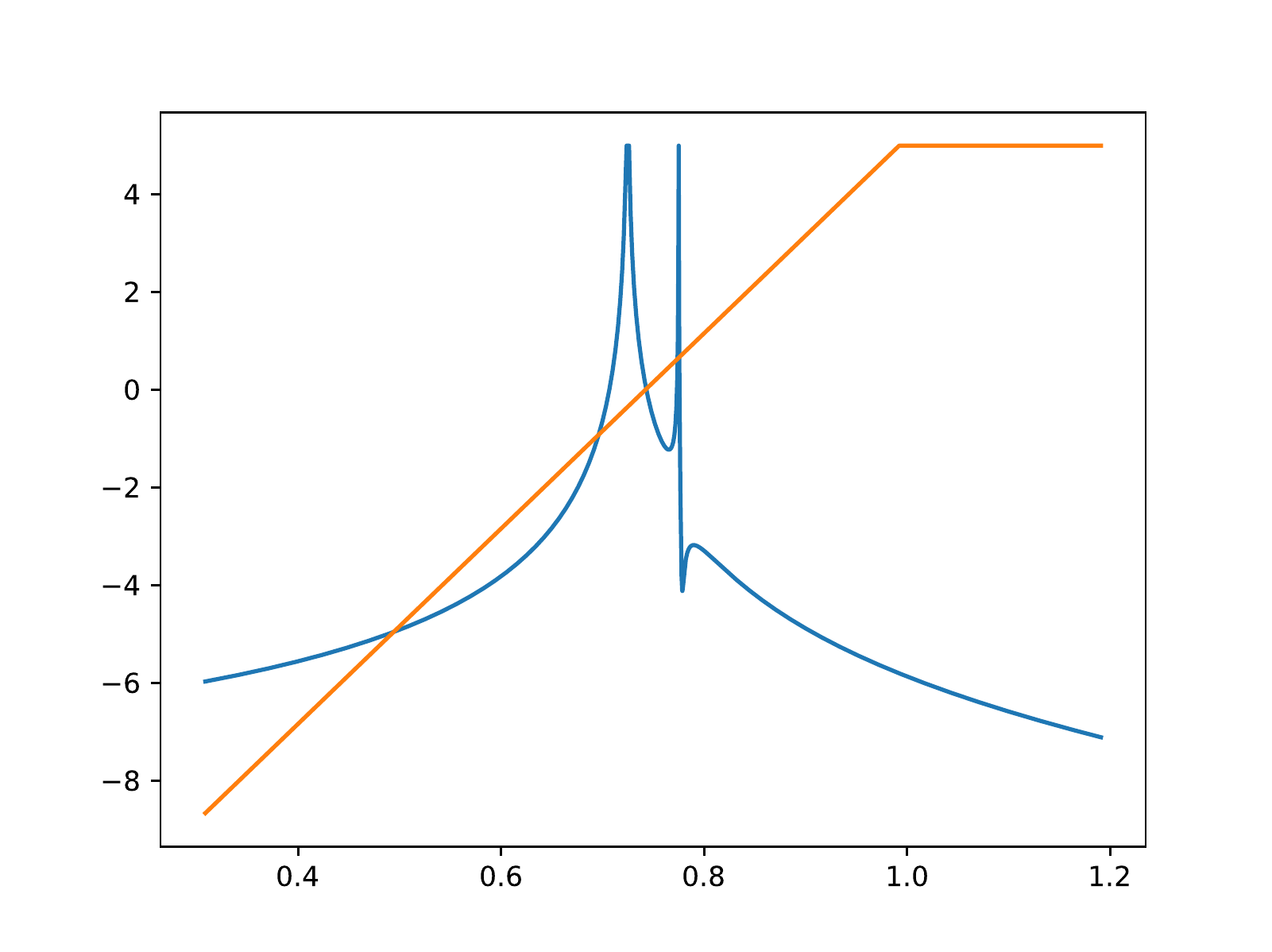}
    \caption{\label{fig:p2.1d2}$S_0$ for $d=2$ and $p=2.1$.}
\end{figure}

\begin{figure}[H]
\centering
 \includegraphics[width=0.3\linewidth]{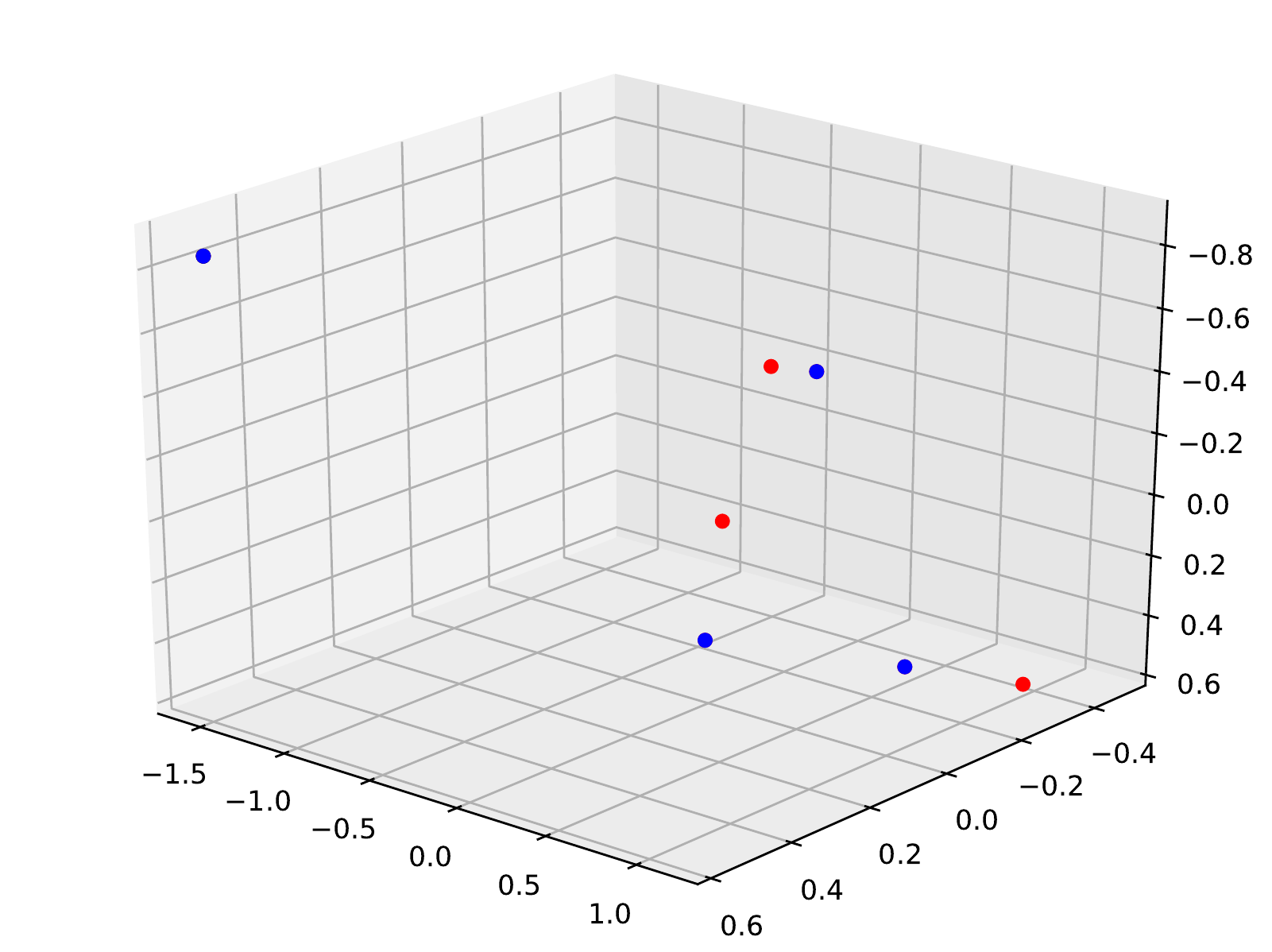}\includegraphics[width=0.3\linewidth]{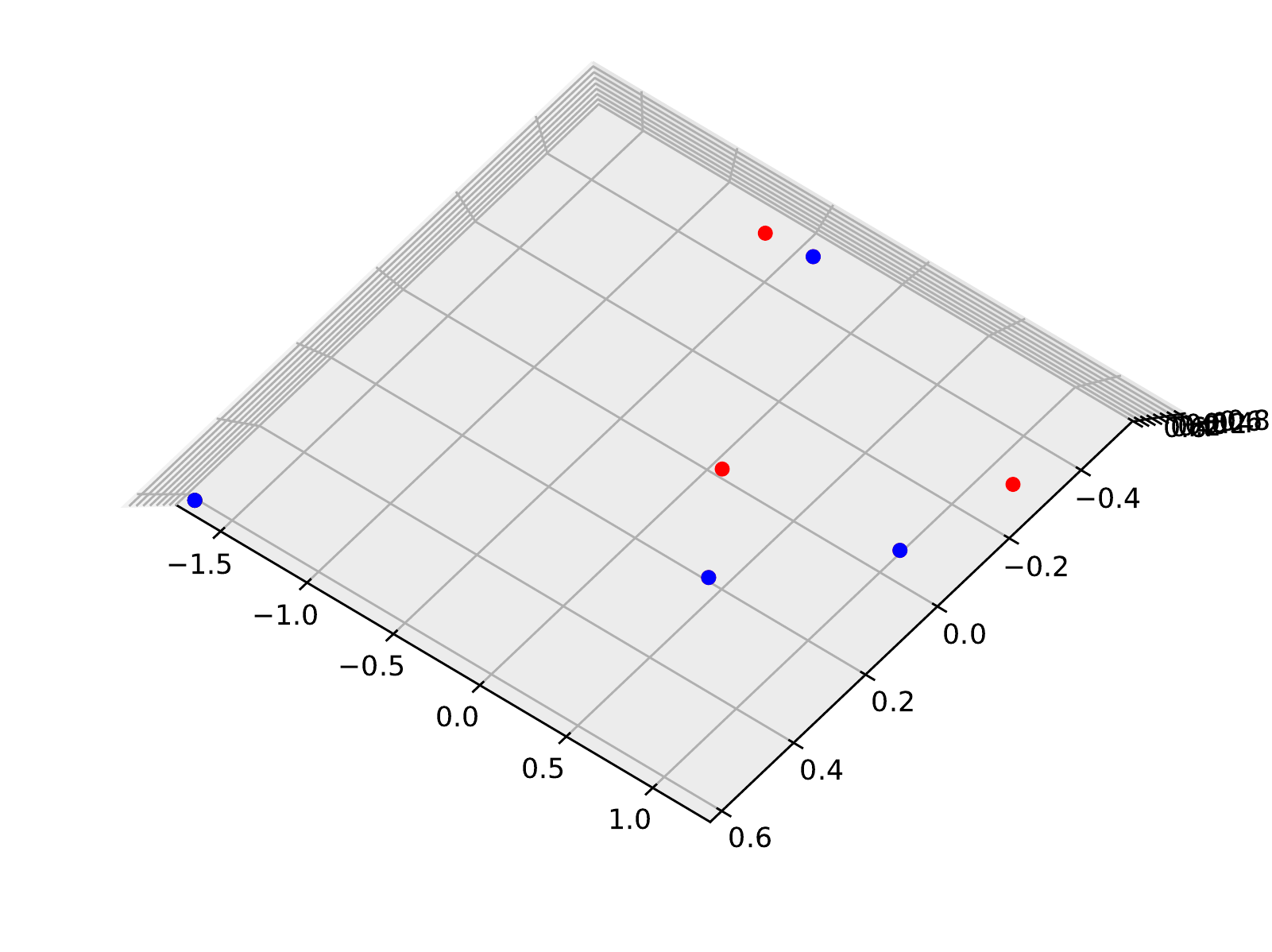}\includegraphics[width=0.3\linewidth]{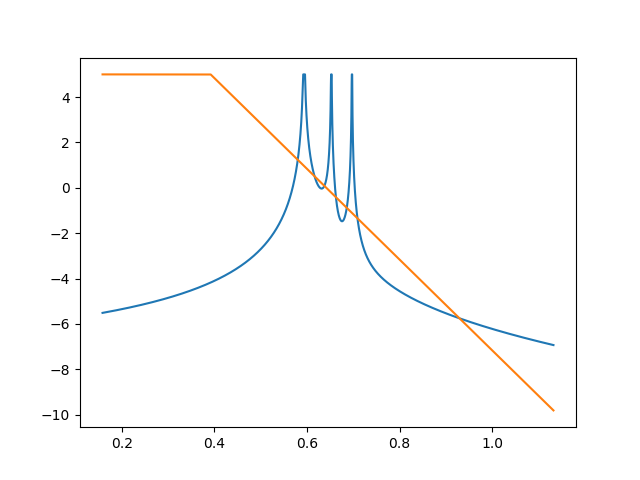}
    \caption{\label{fig:p1.9d3}$S_0$ for $d=3$ and $p=1.9$.}
\end{figure}

\begin{figure}[H]
\centering
 \includegraphics[width=0.3\linewidth]{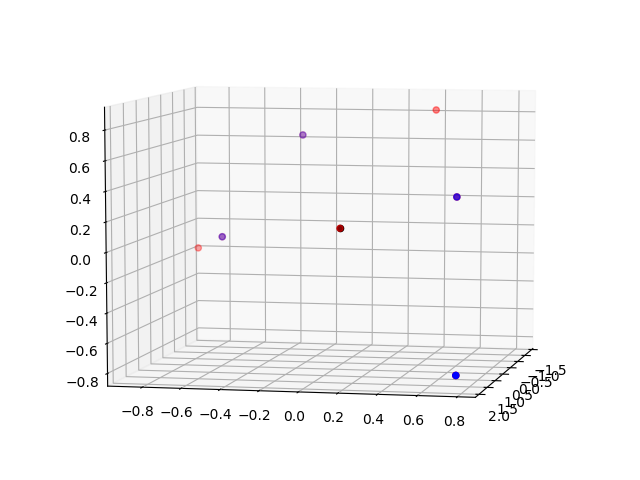}\includegraphics[width=0.3\linewidth]{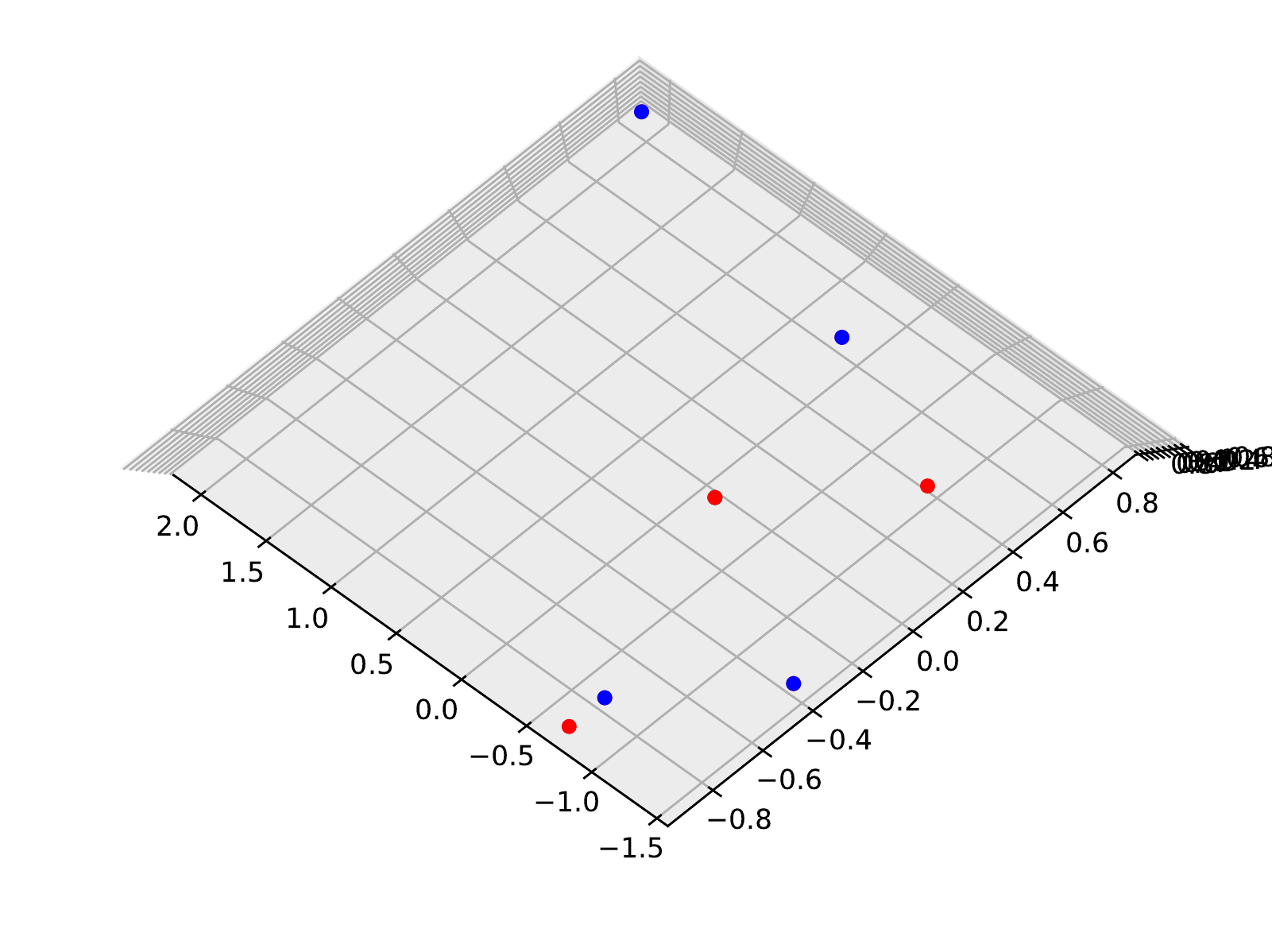}\includegraphics[width=0.3\linewidth]{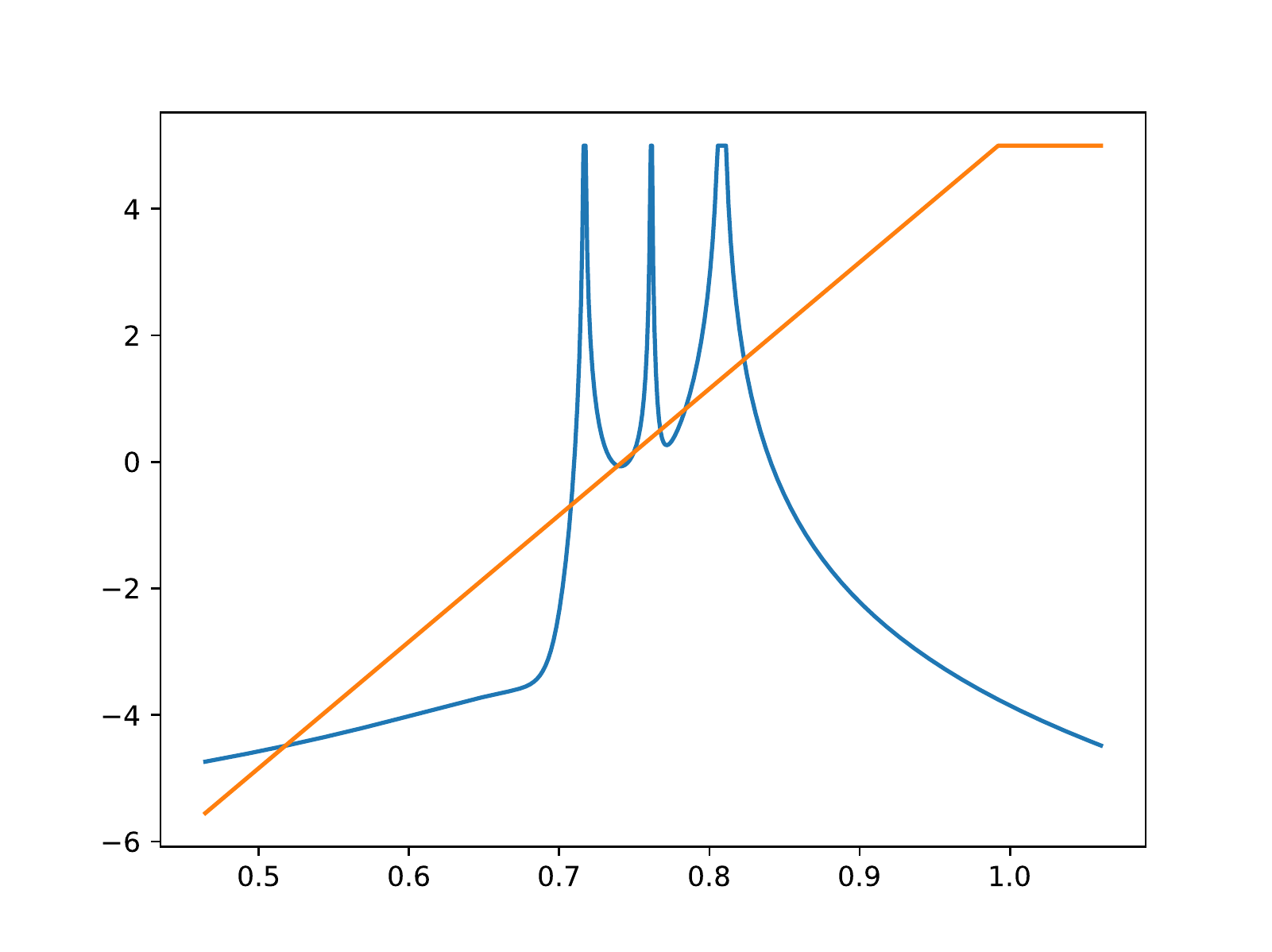}
    \caption{\label{fig:p2.1d3}$S_0$ for $d=3$ and $p=2.1$.}
\end{figure}

\begin{figure}[H]
\centering
 \includegraphics[width=0.45\linewidth]{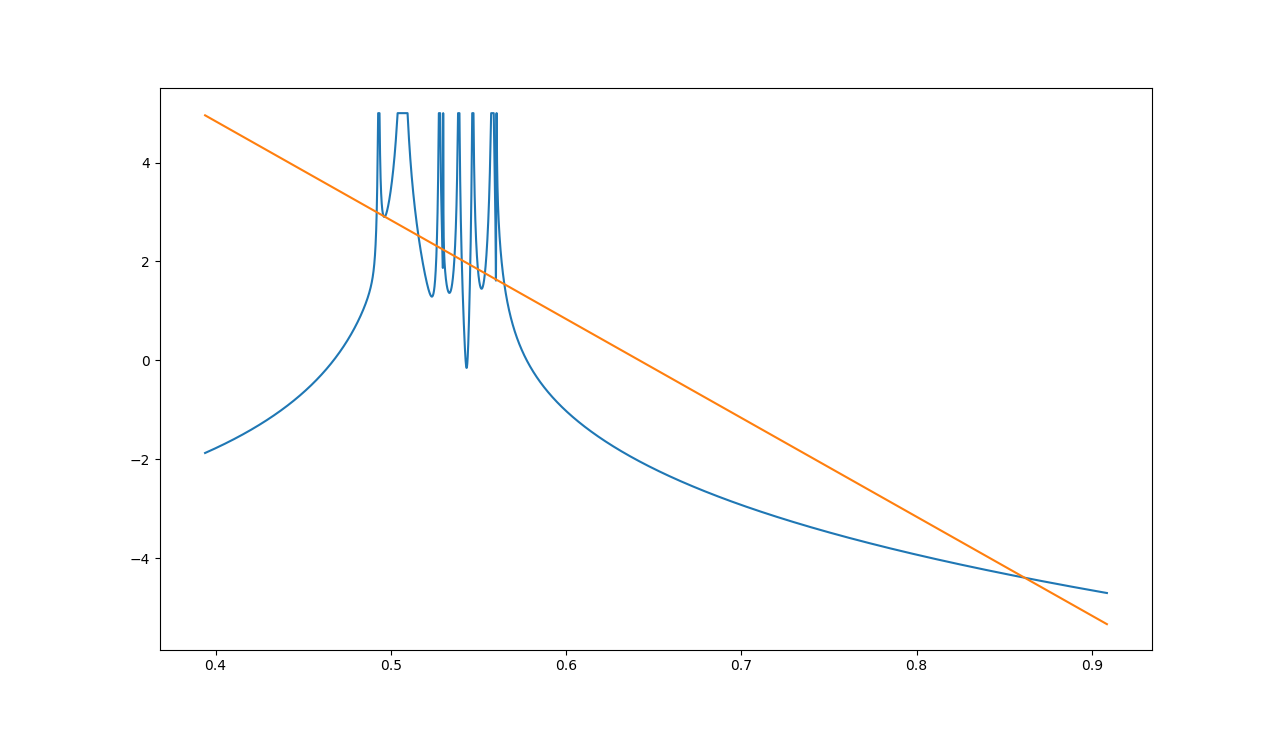}\includegraphics[width=0.45\linewidth]{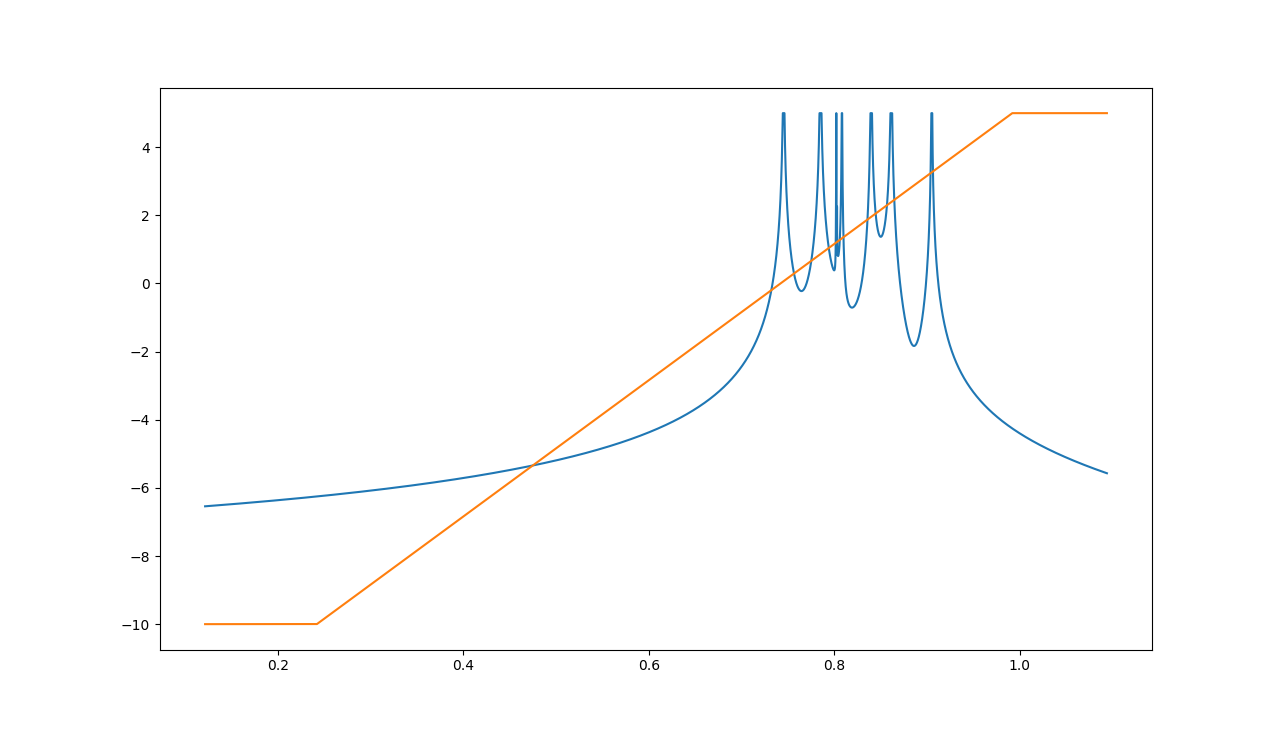}
    \caption{\label{fig:dim8}$S_0$ for $d=8$, $p=1.9$ on the left and $p=2.1$ on the right.}
\end{figure}

In the particular example $c(X,Y)=|X-Y|^p$, the computations are easy as the important unknown parameter $\lambda = p|y|^{p-2}$ is one-dimensional. We coded a solver that generates random $y_1,...,y_{d+1}\in\R^d$ and determines the missing $y_{d+2},...,y_{k}$, with $k=2d$ if $p\le 1$, and $k=2d+1$ if $p>1$ such that $\{y_1,...,y_k\} = \{c_x(0,Y) = A(Y)\}$ for some $A\in \Aff_d$, see Theorem \ref{thm:pdistance}. (As we chose randomly these vectors, we are in a non-degenerate case with probability $1$). Theorem \ref{thm:choquetconv} only covers the case in which $p < 2-\frac25$ or $p>2+\frac23$, however the numerical experiment seems to show that the result of this theorem still holds for all $2\neq p>1$. Figures \ref{fig:p1.9d2}, \ref{fig:p2.1d2}, \ref{fig:p1.9d3}, \ref{fig:p2.1d3}, and \ref{fig:dim8} show configurations ($S_0$, on the left) for $p=1.9$ and $p=2.1$ in which the result of the theorem holds, and the graphs of $\frac{1}{p-2}\log\left(\frac{\lambda}{p}\right)$ compared to $\log\big(\mathbf{y}(-p\lambda^{p-2})\big)$ as functions of $\log(\lambda)$ (on the right). The intersections are in bijection with the points in $S_0$ because of the non-degeneracy by Theorem \ref{thm:fundistance} with the change of variable $t = -p\lambda^{p-2}$. The color of the points need to be interpreted as follows: $d+1$ blue points are chosen at random so that $0$ belongs to their convex hull. Then the new $d$ points given by Theorem \ref{thm:fundistance} are colored in red. Finally the point corresponding to the first intersection of the curves on the right is colored in yellow because this special intersection differentiates the case $p\le 1$ and the case $p>1$. We begin with Figures \ref{fig:p1.9d2} and \ref{fig:p2.1d2}, in two dimensions.

Now Figures \ref{fig:p1.9d3} and \ref{fig:p2.1d3}, in three dimensions.

Finally, Figure \ref{fig:dim8} shows two experiments in which $|S_0|$ contains exactly $17$ elements for $d=8$.

\bibliographystyle{plain}
\bibliography{mabib}
\end{document}